\input amstex
\documentstyle{amsppt}
%\magnification=\magstep2
 \NoRunningHeads
\topmatter
\title
Ergodic cocycles of IDPFT systems~and nonsingular Gaussian actions
\endtitle

\author
Alexandre I. Danilenko and Mariusz Lema{\'n}czyk
\endauthor

%\thanks{Research of the two authors is supported in part by the special program of invitations of the semester
%``Ergodic Theory and Dynamical Systems in their Interactions with Arithmetic and Combinatorics'', Chair Jean Morlet, 1.08.2016-30.01.2017.}
%\endthanks

\thanks{The main results of the paper were obtained during the visit of the first named author at Nicolaus Copernicus University in March-April 2020. The  staying was supported by a special research grant of NCU.}
\endthanks
\thanks
Research of the second named author was supported by Narodowe Centrum Nauki grant UMO-2019/33/B/ST1/00364.
\endthanks

\address
B. I. Verkin Institute for Low Temperature Physics
\& Engineering of National Academy of Sciences of Ukraine, 47 Nauky Ave.,
 Kharkiv, 61103, UKRAINE
\endaddress
\email alexandre.danilenko\@gmail.com
\endemail

\address
Faculty of Mathematics and Computer Science, Nicolaus Copernicus
University, ul. Chopina 12/18, 87-100 Toru\'n, Poland
\endaddress

\email
mlem\@mat.umk.pl
\endemail

\dedicatory 
To the memory of Sergiy Sinel'shchikov, our colleague and friend 
\enddedicatory

\NoBlackBoxes

\abstract
It is proved that each Gaussian cocycle over a mildly mixing Gaussian transformation
is either a Gaussian coboundary or sharply weak mixing.
The class of nonsingular  infinite direct products $T$ of  transformations $T_n$, $n\in\Bbb N$, of finite type (IDPFT) is studied.
It is shown that if $T_n$ is mildly mixing, $n\in\Bbb N$, the sequence of the Radon-Nikodym derivatives of $T_n$ is asymptotically translation quasi-invariant  and $T$ is conservative then the Maharam extension of $T$ is sharply weak mixing. This techniques
provides a new approach to  
the nonsingular Gaussian transformations studied recently by Arano, Isono and Marrakchi.
\endabstract

\endtopmatter

\document

\head 0. Introduction
\endhead
The original motivation of this paper was to tackle a problem (stated in \cite{LeLeSk}) that is related to the theory of Gaussian dynamical systems: let $T$ be an ergodic (equivalently, weakly mixing) Gaussian transformation on a standard probability space $(X,\goth B,\mu)$ and let $H$ be the  corresponding  invariant Gaussian  subspace of the real  Hilbert space $L^2_0(X,\mu)$.

\proclaim{Conjecture}
For each function $f\in H$, either $f$ is a $T$-coboundary (equivalently, a Gaussian coboundary) or the skew product transformation $T_f$ acting on $X\times\Bbb R$ is ergodic.
 \endproclaim

In this paper we obtain the affirmative answer under a slightly stronger assumption than the weak mixing.
We say that a nonsingular transformation $R$ is {\it sharply weak mixing} if the direct product of $R$ with each ergodic conservative transformation is either totally dissipative or ergodic.
In particular, $R$ is ergodic.
We also recall that $T_f$ is conservative for each $f\in H$.

\proclaim{Theorem 0.1} If $T$ is mildly mixing and $f$ is not  a coboundary for $T$ then $T_f$ is sharply weak mixing.
\endproclaim

To prove Theorem~0.1 we note that there exists a decomposition of $T$ into direct product
of mildly mixing transformations $T_n$ in a such a way that
 $f$ splits into
a sum of coboundaries $f_n:=a_n-a_n\circ T_n$ for $T_n$, $n\in\Bbb N$. 
Moreover, the sequence of distributions of the transfer functions $(a_n)_{n=1}^\infty$ satisfies a certain property that we call  {\it ATI (asymptotic translation invariance)} in~Definition~1.2.
Then Theorem~0.1 follows from  the next theorem.

\proclaim{Theorem 0.2} Given a locally compact second countable Abelian group $G$,  a sequence of mildly mixing dynamical systems $(X_n,\nu_n, T_n)$   and a sequence of functions $f_n:X_n\to G$, $n\in\Bbb N$, consider the infinite direct product $(X,\nu,T):=\bigotimes_{n=1}^\infty (X_n,\nu_n,T_n)$.
Suppose that a function $f(x):=\sum_{n=1}^\infty(f_n(T_nx_x)-f_n(x_n))\in G$ is well defined for $\nu$-a.e. $x=(x_n)_{n=1}^\infty\in X$.
If the sequence of distributions $(\nu_n\circ f_n^{-1})_{n\in\Bbb N}$ on $G$ is ATI and the $f$-skew product extension $T_f:X\times G\to X\times G$ of $T$ is conservative then $T_f$ is sharply weak  mixing.
\endproclaim

The proof of Theorem~0.2 is based on the two ideas: 
\roster
\item"---" the mild mixing and the  product structure of $T_f$  yield  that each $T_f$-invariant subset is also invariant under a large group of  ``finitary'' transformations, i.e. transformations that ``move'' finitely many coordinates only;
\item"---" the ATI property implies that this finitary group is ergodic via techniques related to computation of the essential values of  cocycles.
\endroster
The first aforementioned idea was inspired by  the proof \cite{ArIsMa, Theorem~D} on
ergodic properties of some nonsingular Gaussian group actions.

We then turn to classical problems of nonsingular ergodic theory.
We mention a recent progress in providing natural examples for nonsingular ergodic theory:  nonsingular Bernoulli and Markov shiftwise actions (see  \cite{DaLe}, \cite{KoSo}, \cite{Av}, \cite{MaVa}
and references therein), nonsingular Gaussian systems \cite{ArIsMa}, nonsingular Poisson systems (\cite{DaKoRo1}, \cite{DaKoRo2}).
In the present work we introduce one more natural family of nonsingular transformations. 
We say that  a nonsingular transformation $T$ on a standard probability space $(X,\mu)$ is an {\it infinite direct product  of finite types (IDPFT)} if there is 
a sequence of ergodic probability preserving  dynamical systems $(X_n,\nu_n,T_n)$ and a sequence of probability measures
$\mu_n$ on $X_n$, $n\in\Bbb N$, such that $\mu_n\sim\nu_n$ for each $n$
and $(X,\mu,T)=\bigotimes_{n=1}^\infty(X_n,\mu_n,T_n)$.
Kakutani's theorem \cite{Ka} provides a criterion where $\mu$ is quasi-invariant under $T$.
We are interested in the case where $\mu\perp\nu$
and $\mu$ does not admit an equivalent  $T$-invariant probability.
It is possible that $(X,\mu, T)$ is totally dissipative.
Moreover, we show that for each ergodic conservative nonsingular transformation $S$, the product $T\times S$ is either  totally dissipative or conservative.

\proclaim{Theorem 0.3} Let $(X_n,\nu_n,T_n)$ be mildly mixing for each $n>0$.
If $T$ is $\mu$-conservative and the sequence of distributions of the random variables $\log\frac{d\mu_n}{d\nu_n}$, $n\in\Bbb N$, is ATQI then $T$ is ergodic of stable type Krieger's type $III_1$.
Moreover, the Maharam extension of $T$ is sharply weak mixing.
\endproclaim

The property ATQI (asymptotically translation quasi-invariantness, see Definition~2.8) in the statement of Theorem~0.3 is an analogue of ATI though neither ATI implies ATQI nor vice-versa. 
The scheme of the proof of Theorem~0.3 is similar to the that of Theorem~0.2 and we
use again the aforementioned two ideas. 
However, there is a ``nonsingular'' nuance. 
Namely, 
a formal repetition of  the proof of Theorem~0.2 yields that  the  group of finitary transformations is ergodic with respect to the ``wrong'' measure.
Hence, it does not work.
We recall that there are two different (mutually singular) natural measures associated with a IDPFT system:   $\nu$ (invariant) and  $\mu$ (quasiinvariant). 
Therefore a certain additional argument  and the property ATQI instead of ATI  are needed
to prove ergodicity for the ``right'' measure.
We also provide examples of rigid IDPFT systems $T$ of Krieger's type $III_\lambda$ for an arbitrary $\lambda\in(0,1)$.

We have already mentioned that the nonsingular Gaussian systems were studied recently in \cite{ArIsMa}.
However, the exposition there is based heavily  on the affine geometry and  often use a nonstandard (from the dynamical viewpoint) terminology. 
Therefore, we decided to provide here an alternative exposition of this important topic.
We define the nonsingular Gaussian systems as transformations on Hilbert spaces $\Cal H$ furnished with Gaussian measures 
stressing on the fact  that
the systems  are compositions of classical Guassian automorphisms and totally dissipative transformations (given by nonsingular rotations).
Connections with the underlying Fock space, the first chaos and  the exponential map
are enlighten explicitly.
We also explain interrelation between the nonsingular Gaussian systems and the nonsingular Poisson systems.
Our main observation is that the  Gaussian transformations (out of a ``small'' family of degenerated ones) is a subclass of  IDPFT systems.
Hence we deduce from Theorem~0.2 one of the main results of \cite{ArIsMa}
($\Cal H_0$ below is a linear subspace of $\Cal H$ endowed with a new inner product, see Section~3).

\proclaim{Theorem 0.4} Let  an orthogonal operator $V$ of a real Hilbert space $\Cal H_0$
be mildly mixing. 
Let $f\in\Cal H_0$ not be  a $V$-coboundary, i.e. $f\ne Va-a$ for any $a\in\Cal H_0$.
If the nonsingular Gaussian transformation $T_{(f,V)}$ associated with the pair $(f,V)$
is conservative then the Maharam extension of $T_{(f,V)}$ is sharply weak mixing. 
In particular, $T_{(f,V)}$ is of type $III_1$.
\endproclaim

The outline of the paper is as follows. In Section~1 we introduce main definitions: Hellinger distance, weak mixing properties for nonsingular actions,  ATI property, skew product extension, essential value of a cocycle, etc. 
Then we prove Theorem~0.2 (see Theorem~1.5) and deduce Theorem~0.1 from it (see Theorem~1.6).
We also provide a generalization of Theorem~0.1 (see Conjecture~II and a discussion above it).
In Section~2 we consider nonsingular versions of the problems studied in \S1.
IDPFT systems are introduced in Definition~2.2.
Radon-Nikodym cocycle, Maharam extension and Kriger's types $III_\lambda$, $0\le\lambda\le 1$, are discussed there.
We show that each IDPFT system is either conservative or totally dissipative (Corollary~2.7),
introduce the ATQI property (Definition~2.8) and prove Theorem~0.3 (Theorem~2.10).
Type $III_\lambda$ rigid IDPFT systems are also constructed there for each $\lambda\in(0.1)$ (Proposition~2.12).
The final Section~3 is devoted to nonsingular Gaussian systems.
We first recall the definition of Gaussian measure in a separable Hilbert space.
Then we discuss the main properties of the related Fock space and exponential map.
Given an orthogonal operator $V$ in a Hilbert space $\Cal H_0$ and a vector $f\in\Cal H_0$,
we associate a nonsingular transformation $T_{(f,V)}$ acting on the corresponding Hilbert space $\Cal H\supset \Cal H_0$ equipped with a Gaussian measure $\mu$. 
We show that $T_{(f,V)}$ is the composition of the classic Gaussian $\mu$-preserving transformation associated to $V$ with the (totally dissipative) rotation by $f$.
It is well known that the nonsingular transformation  group $\{T_{(f,0)}\mid f\in\Cal H_0\}$
generated by the rotations is ergodic (see, e.g. \cite{Gu}) but the Kriger's type  has  not been specified so far. 
We prove that it is $III_1$ (Theorem~3.7).
We show that the Koopman operator generated by   $T_{(f,V)}$ is the Weyl operator associated to the pair $(f/2,V)$.
A criterion for the existence of an invariant equivalent probability measure for  $T_{(f,V)}$ is established
 in Theorem~3.9 (cf.  \cite{DAKoRo1, Proposition~6.4} and \cite{ArIsMa}).
 Theorem~0.4 in proved in this section (Theorem~3.12).

 After completion of this paper we learnt about a work \cite{MaVa}  devoted to nonsingular Gaussian actions of arbitrary groups. It  was written independently but simultaneously with our work.\footnote{The two papers appeared on ArXiv in two successive days. We thank S.~Vaes for informing us about \cite{MaVa}.} Some of our  results overlap with theirs.
 Say, Theorem~3.7 is  \cite{MaVa, Theorem~3.1} and
 Theorem~0.1, though stated in a more general form, is, in fact, equivalent to \cite{MaVa, Theorem~9.1(3)} in case of $\Bbb Z$-actions.
 Our proofs are different.
 They are based solely on elementary techniques of the nonsingular ergodic (measurable orbit) theory.
 We do not use  affine geometry, representation theory nor harmonic analysis.

\head 1. Weak mixing cocycles of product type.
\endhead
\subhead{1.1. Hellinger distance and Kakutani's theorem}
\endsubhead
Let $\gamma$ and $\delta$ be two equivalent  probability measures on a standard Borel space $(Y,\goth C)$.
The square of the Hellinger distance  between $\gamma$ and $\delta$ is
$$
H^2(\gamma,\delta):=\frac12\int_Y\bigg(1-\sqrt{\frac{d\gamma}{d\delta}}\bigg)^2\,d\delta=1-\int_Y\sqrt{\frac{d\gamma}{d\delta}}d\delta.
$$
By the Cauchy--Schwarz inequality, $0\le H(\gamma,\delta)<1$.
We also remind \cite{Ni} the following inequalities between the Hellinger distance and the total variation:
$$
H^2(\gamma,\delta)\le\|\gamma-\delta\|_1:=\sup_{C\in\goth C}|\gamma(C)-\delta(C)|\le \sqrt 2H(\gamma,\delta).\tag1-1
$$
We now state the Kakutani theorem on equivalence of infinite products of probability measures \cite{Ka}.

\proclaim{Theorem A} Let $\mu_n$ and $\nu_n$ be two equivalent probability measures on a standard Borel space $(X_n,\goth B_n)$ for each $n\in\Bbb N$.
Let $\mu$ and $\nu$ denote the infinite product measures $\bigotimes_{n\in\Bbb N}\mu_n$
and $\bigotimes_{n\in\Bbb N}\nu_n$ respectively on the standard Borel space $(X,\goth B):=\bigotimes_{n\in\Bbb N}(X_n,\goth B_n)$.
If 
$$
\prod_{n=1}^\infty \Big(1- H^2(\mu_n,\nu_n)\Big)>0\text{ \ or, equivalently, \  }\sum_{n=1}^\infty H^2(\mu_n,\nu_n)<\infty\tag1-2
$$
then $\mu\sim\nu$, $\prod_{n=1}^\infty (1- H^2(\mu_n,\nu_n))=1-H^2(\mu,\nu)$ and
$
\frac{d\mu}{d\nu}(x)=\prod_{n\in\Bbb N}\frac{d\mu_n}{d\nu_n}(x_n)
$
at a.e. $x=(x_n)_{n\in\Bbb N}\in X$.
If \thetag{1-2} does not hold then $\mu\perp\nu$.
\endproclaim

\subhead 1.2. Weak mixing properties of nonsingular actions
\endsubhead
We recall that given a non-singular transformation $R$ of a standard Borel probability space $(Y,\goth C,\nu)$, there is a unique decomposition  $Y=\Cal D(R)\sqcup \Cal C(R)$ (called {\it Hopf's decomposition}) of $Y$ into two Borel sets such that $\Cal D(R)$ is the disjoint union of the orbit of a {\it wandering set} $W$, i.e. $\Cal D(R) =\bigsqcup_{n\in\Bbb Z} R^nW$  and $\Cal C(R)=Y\setminus \Cal D(R)$ contains no non-trivial wandering set. 
If $\Cal C(R)=Y$ then $R$ is called {\it conservative} and if $\Cal D(R)=Y$ then $R$ is called {\it totally dissipative}.
 As both parts $\Cal C(R)$ and $\Cal D(R)$ are $R$-invariant, each ergodic $R$ is either conservative or totally dissipative. 
An ergodic  conservative nonsingular transformation $R$ is called {\it weakly mixing} if for each ergodic probability preserving transformation $S$, the Cartesian product $R\times S$ 
is ergodic.
We now introduce a stronger concept of weak mixing.

\definition{Definition 1.1}
An ergodic conservative nonsingular transformation $R$ is called {\it sharply weak mixing} if for each ergodic conservative nonsingular transformation $S$, the direct product $R\times S$ is either
totally dissipative or ergodic.
\enddefinition

If $S$ in the above definition admits an equivalent invariant probability measure (i.e. $S$ is of type $II_1$) then $T \times S$ is conservative (see   \cite{Aa, Proposition~1.1.6, part 2}).
Hence $R\times S$ is ergodic according to Definition~1.1.
Thus,
 every sharply weak mixing transformation is weakly mixing.
 It follows from \cite{SiTh} that every conservative nonsingular transformation with property $K$
is sharply weak mixing (see also \cite{AaLiWe, Theorem~6.7} for other examples).
In \cite{AdFrSi} and \cite{Da} examples of  weakly mixing infinite measure preserving  rank-one transformations $R$
were constructed such that $R\times R$ is conservative but not ergodic.
Hence $R$ is not sharply weak mixing.
We recall that an ergodic probability preserving transformation $R$ defined on a space $(Y,\goth C,\nu)$ is called {\it mildly mixing} (\cite{FuWe}, see also \cite{AaLiWe} and \cite{ScWa}) if every function $f\in L^\infty(\nu)$ such that $\|f\circ T^{n_i}-f\|_1\to 0$ for some sequence $n_i\to\infty$  is constant.

 We will utilize the following result from \cite{ScWa}.
 
 \proclaim{Theorem B}
 Let $R$ be a mildly  mixing transformation of a standard probability
 space $(Y,\goth C,\nu)$ and let $C$ be a conservative nonsingular transformation of a standard probability space  $(Z,\goth F,\tau)$.
 If a function $F\in L^\infty(Y\times Z,\nu\otimes\tau)$
 is invariant under $R\times C$ then there is $f\in L^\infty(Z,\tau)$ such that $F(y,z)=f(z)$ a.e.
 \endproclaim

We note that Theorem~B was proved in \cite{ScWa} for the ergodic conservative $C$ only but the proof remains valid for an arbitrary conservative $C$ as well.
Direct products of finitely (and countably) many mildly mixing transformatations are mildly mixing.

It follows from Theorem~B that  an ergodic finite measure preserving transformation is
 sharply weak mixing if and only if it is mildly mixing.
In Theorems~1.5 and 1.6 below we will provide examples of mildly mixing transformations (including zero entropy case) which have  locally compact group extensions that are sharply weak mixing infinite measure preserving (and hence not mildly mixing).

\subhead 1.3. ATI property
\endsubhead
Fix a locally compact second countable Abelian group $G$.
Denote by $\lambda_G$ a Haar measure on $G$.

\definition{Definition 1.2}
A sequence $(\xi_n)_{n=1}^\infty$ of probability Borel measures on $G$ is called
{\it asymptotically translation invariant} (ATI) if 
$$
\lim_{m\to\infty}\|\xi_n*\xi_{n+1}*\cdots*\xi_{n+m}*\delta_a-\xi_n*\xi_{n+1}*\cdots*\xi_{n+m}\|_1= 0
$$ 
  for each $n\in\Bbb N$ and $a\in G$.
\enddefinition

\example{Example 1.3}
Let $\Cal N_{a,\sigma^2}$ denote the normal distribution on $\Bbb R$ with parameters $a$ and $\sigma^2$, i.e.
$\widehat{\Cal N_{a,\sigma^2}}(t)=e^{iat-\frac12\sigma^2t^2}$ for all $t\in\Bbb R$.
Given two sequences $(a_n)_{n\in\Bbb Z}$ and $(\sigma_n)_{n=1}^\infty$ of reals such that $\sum_{n=1}^\infty\sigma_n^2=+\infty$, the sequence of probabilities
$(\Cal N_{a_n,\sigma_n^2})_{n=1}^\infty$ is ATI.
Indeed, 
$$
\align
\Cal N_{a_n,\sigma_n^2}*\cdots *\Cal N_{a_{n+m},\sigma_{n+m}^2}&=
\Cal N_{\sum_{k=n}^{n+m}a_k,{\sum_{k=n}^{n+m}\sigma^2_k}},\\
\Cal N_{a_n,\sigma_n^2}*\cdots *\Cal N_{a_{n+m},\sigma_{n+m}^2}*\delta_a&=
\Cal N_{a+\sum_{k=n}^{n+m}a_k,{\sum_{k=n}^{n+m}\sigma^2_k}}\quad\text{and}\\
H^2\Big(\Cal N_{a+\sum_{k=n}^{n+m}a_k,\sum_{k=n}^{n+m}\sigma^2_k},\Cal N_{\sum_{k=n}^{n+m}a_k,\sum_{k=n}^{n+m}\sigma^2_k}\Big)&=1-e^{-\frac18\cdot\frac{a^2}{\sum_{k=n}^{n+m}\sigma^2_k}}\to 0
\endalign
$$
as $m\to\infty$\footnote{We use the fact that $H^2(\Cal N_{a,\sigma^2},\Cal N_{b,\tau^2})=1-
\sqrt{\frac{2\sigma\tau}{\sigma^2+\tau^2}}\,e^{-\frac14
\frac{(a-b)^2}{\sigma^2+\tau^2}}$ for all $a,b,\sigma,\tau\in\Bbb R$.}.
Hence $(\xi_n)_{n=1}^\infty$  is ATI in view of  \thetag{1-1}.
\endexample

\subhead 1.4. Ergodic cocycles of ergodic transformation groups
\endsubhead Given a standard Borel $\sigma$-finite measure space $(Y,\goth C,\nu)$, we denote by Aut$(Y,\nu)$ the group of all $\nu$-nonsingular invertible Borel transfomations on $Y$.
Let Aut$_0(Y,\nu)$ denote the subgroup of $\nu$-preserving transformations from
Aut$(Y,\nu)$.
Let $\Gamma$ be an ergodic countable subgroup in Aut$(Y,\nu)$.
The  {\it full group} $[\Gamma]$ of $\Gamma$ is defined by:
$$
[\Gamma]:=\{\theta\in\text{Aut}(Y,\nu)\mid \theta y\in\{\gamma y\mid\gamma\in \Gamma\}\text{ at a.e. }y\in Y\}.
$$
A measurable map $\alpha:\Gamma\times Y\to G$ is called a {\it cocycle} of $\Gamma$ if
$$
\alpha(\gamma_1\gamma_2,y)=\alpha(\gamma_1,\gamma_2y)+\alpha(\gamma_2,y)\quad\text{ at a.e. }y\in Y\tag1-3
$$
for all $\gamma_1,\gamma_2\in\Gamma$.
From now on  we assume that $\Gamma$ is free, i.e. if $\gamma\in\Gamma\setminus\{I\}$
then $\gamma y \ne y$ for a.e. $y$.
Then $\alpha$ can be ``extended'' to $[\Gamma]$ if we set
$$
\alpha(\theta,y):=\alpha(\gamma,y)\quad\text{where $\gamma$ is defined uniquely by }\theta y=\gamma y.
$$
It is straightforward to verify that \thetag{1-3} holds if we replace $\gamma_1$ and $\gamma_2$ with arbitrary elements from $[\Gamma]$.
A cocycle $\alpha$ is a {\it coboundary} if there is 
a measurable map $a:Y\to G$ such that
$$
\alpha(\gamma,y)=a(\gamma y)-a(y)\quad\text{ at a.e. }y\in Y
$$
for all $\gamma\in \Gamma$.
Given a pair $(\Gamma,\alpha)$, we can construct a transformation group $\Gamma_\alpha:=\{\gamma_\alpha\mid\gamma\in\Gamma\}\subset\text{Aut}(Y\times G,\nu\times\lambda_G)$, where
$$
\gamma_\alpha(y,g):=(\gamma y,\alpha(\gamma,y)+g)\quad\text{ for all }y\in Y, g\in G.
$$
The group $\Gamma_\alpha$ is called the {\it $\alpha$-skew product extension} of $\Gamma$.
If $\Gamma$ preserves $\nu$ then $\Gamma_\alpha$  preserves the product measure $\nu\otimes\lambda_G$.
If $\Gamma_\alpha$ is ergodic then $\alpha$ is called {\it ergodic}.
A coboundary is never ergodic (unless $G$ is a singleton).
It is easy to verify that if  $\Gamma=\{R^n\mid n\in\Bbb Z\}$ for a transformation $R\in\text{Aut}(Y,\nu)$ then each measurable function $f:Y\to G$ defines uniquely  a cocycle $\alpha_f$ of $\Gamma$ via the condition
$$
\alpha_f(R,y):=f(y)\quad\text{for each }y\in Y.
$$
For brevity we will write $R_f$ for the $\alpha_f$-skew product extension $R_{\alpha_f}$ of $R$.

We now recall an important concept of  essential value for  a cocycle.

\definition{Definition 1.4} Suppose that $\Gamma$ preserves $\nu$.
An element $g\in G$ is called an {\it essential value of $\alpha$}
if for each subset $A\subset Y$ of positive measure and a neighborhood $U$ of $g$, there are a Borel subset $B\subset A$ and an element $\gamma\in\Gamma$ such that $\nu(B)>0$,
$\gamma B\subset A$ and $\alpha(\gamma, y)\in U$ for all $y\in B$.
\enddefinition

It appears that the set $r(\alpha)$ of all essential values of a cocycle is a closed subgroup of $G$.
Our interest to the essential values of $\alpha$ is explained by the fact that
 $\alpha$ is ergodic if and only if $r(\alpha)=G$ \cite{Sc}.
It is often easier to check the aforementioned condition on essential values  not for each subset $A\in\goth C$ of positive measure but only  for a dense subfamily of subsets in $\goth C$.
However in this case we have to strengthen this condition.
More precisely, 
we will use the
 following  lemma.
(It follows, for example, from \cite{Da, Lemma~2.4}.)

\proclaim{Lemma C} Let $(Y,\goth C,\nu)$ be a standard probability space,  $\goth A$ a dense subset in $\goth C$, $\Gamma$  an ergodic  countable subgroup of $\text{\rom{Aut}}_0(Y,\nu)$ and
$\alpha:\Gamma\times Y\to G$ a Borel cocycle of  $\Gamma$.
If for some $a\in G$ and each subset $B\in\goth A$ and each neighborhood $U$ of $0$ in $G$, there are a measurable subset
$D\subset B$ and an element $\theta\in[\Gamma]$ such that $\theta D\subset B$,
$\nu(D)>0.5\nu(B)$ and $\alpha(\theta,x)\in a+U$ for all $x\in D$ then $a$ is an essential value of $\alpha$.
\endproclaim

\subhead 1.5. Sharp weak mixing of skew products for  cocycles of product type
\endsubhead
 In this subsection we prove the following theorem.

\proclaim{Theorem 1.5} 
Let $T_n$ be a mildly mixing    transformation of a standard probability space $(X_n,\goth B_n,\nu_n)$ for each $n\in\Bbb N$.
Let 
$$(X,\goth B,\nu, T):=\bigotimes_{n\in\Bbb Z}(X_n,\goth B_n,\nu_n,T).
$$
Suppose that for a measurable function $f:X\to G$, there are functions
$f_n:X_n\to G$ such that
 $f(x)=\sum_{n\in\Bbb N}(f_n(T_nx_n)-f_n(x_n))$  at $\nu$-a.e. $x=(x_n)_{n=1}^\infty\in X$
and 
the sequence of measures $(\nu_n\circ f_n^{-1})_{n\in\Bbb N}$ is ATI.
If
the skew product extension  $T_f:X\times G\to X\times G$ of $T$ is conservative then $T_f$ is sharply weak mixing.
\endproclaim

\demo{Proof} Let $C$ be  an ergodic conservative transformation of a standard probability space $(Z,\goth Z,\kappa)$.
Suppose that $(\mu\otimes\lambda_G\otimes\kappa)(\Cal C(T_f\times C))>0$.
Since $T_f\times C$ commutes with $I\times C$,
$$
(I\times C)\Cal C(T_f\times C)=\Cal C(T_f\times C).
$$
Since
 $C$ is ergodic, $\Cal C(T_f\times C)=\Cal C(T_f)\times Z=X\times G\times Z$,
where the latter equality holds because  $T_f$ is conservative.
Thus, 
 the direct product $T_f\times C$ is conservative.
 It remains to show that $T_f\times C$ is ergodic.

  Let a function $F\in L^\infty(X\times G\times Z,\mu\otimes\lambda_G\otimes\kappa)$ be invariant under $T_f\times C$.
  We first show that $F$ is also invariant under a huge group of transformations.
Fix $n>0$.
For each $x\in X$, we write $x_1^n:=(x_1,\dots, x_n)\in X_1\times\cdots \times X_n$ and $x_{n+1}^\infty:=(x_{n+1},x_{n+2}, \dots)\in X_{n+1}\times X_{n+2}\times\cdots$.
Then $x=(x_1^n,x_{n+1}^\infty)$.
We define a    measure preserving automorphism  $E_n$  of $(X\times G\times Z,\mu\otimes\lambda_G\otimes\kappa)$
and a nonsingular automorphism $V_n$
  of $(\bigotimes_{k=n+1}^\infty(X_k,\nu_k))\otimes (G,\lambda_G)\otimes (Z,\kappa )$
respectively
by setting
$$
\align
E_n(x,g,z)&:=\Big(x,g+\sum_{k=1}^nf_k(x_k),z\Big)\quad\text{and}\\
V_n((x_k)_{k=n+1}^\infty,g,z)&:=\Big((T_{k}x_k)_{k=n+1}^\infty, g+\sum_{k>n}(f_k(T_kx_k)-f_k(x_k)),Cz\Big).
\endalign
$$
A straightforward verification shows that
$$
E_n (T_f\times C)E_n^{-1}=(T_1\times\cdots \times T_n)\times V_n.
$$
Since $V_n$ is a factor of the transformation $E_n(T_f\times C)E_n^{-1}$ and the latter transformation is conservative, it follows that $V_n$ is conservative.
On the other hand,  the function $F\circ E_n$ is invariant under 
$E_n^{-1}(T_f\times C)E_n$.
Utilizing these two facts
we deduce  from Theorem~B
that $F\circ E_n$ does not depend on the coordinates $x_1,\dots,x_n$.
Hence, for each transformation $S\in\text{Aut}_0(X_1\times\cdots\times X_n,\bigotimes_{k=1}^n\nu_k)$,
we have that
 $
 F\circ E_n\circ (S\times I)=F\circ E_n.
 $
 Therefore
$F$ is invariant under the transformation $E_n(S\times I) E_n^{-1}\in\text{Aut}_0(X\times G\times Z,\mu\times\lambda_G\times\kappa)$ and
$$
E_n(S\times I) E_n^{-1}(x,g,z)=
(S x_1^n,x_{n+1}^\infty, g-A_n(x^n_1)+A_n(S x_{1}^n), z),\tag1-4
$$
where $A_n$ stands for the mapping $X_1\times\cdots\times X_n\ni(x_1,\dots,x_n)\mapsto \sum_{k=1}^nf_k(x_k)$.
Thus, we have shown that $F$ is invariant under each transformation from the  set
$$
\Cal G:=\bigcup_{n>0}E_n\bigg(\text{Aut}_0\bigg(X_1\times\cdots\times X_n,\bigotimes_{k=1}^n\nu_k\bigg)\times \{I\}\bigg) E_n^{-1}.%\subset\text{Aut}_0(X,\nu).
$$
We now consider a new dynamical system.
The space of this system is the product $(X,\goth B,\mu)$.
Denote by $\Gamma$ the group of  transformations of this space generated by mutually commuting measure preserving transformations $\widehat T_1,\widehat T_2,\dots$, %and $\widehat C$, 
where
$$
\align
\widehat T_nx&=(x_{1}^{n-1}, T_nx_n,x_{n+1}^{\infty}),\quad n\in\Bbb N.\\
%\widehat C((x_i)_{i=1}^{\infty},z)&=((x_i)_{i=1}^{\infty}, Cz).
\endalign
$$
Then $\Gamma$ is  countable, Abelian\footnote{It is isomorphic to $\bigoplus_{n=1}^\infty\Bbb Z$.}  and ergodic.
%Denote by $\Cal R$ the $\Gamma$-orbit equivalence relation.
%We define a cocycle $\alpha$ of $\Cal R$ to $G$ by setting:
For each $n>0$, we consider a  coboundary
$$
\alpha_n:X\ni x\mapsto \alpha_n(x):=f_n(T_nx_n)-f_n(x_n)\in G
$$
of $\widehat T_n$.
It is straightforward to verify\footnote{This follows from the fact that each function $\alpha_n$ depends only on a single coordinate $x_n$, $n=1,2,\dots.$} that the $\alpha_n$-skew product extensions
$(\widehat T_n)_{\alpha_n}$ of $\widehat T_n$, $n\in\Bbb N$,   commute mutually.
It follows that a cocycle $\alpha:\Gamma\times X\to G$ of $\Gamma$ with values in $G$ is well defined by the following formulae: 
$$
\alpha(\widehat T_n,x):=\alpha_n(x),\quad n\in\Bbb N.
$$
%However $\alpha$ is not a coboundary.
Since $\alpha_n(x)=A_n((I\times T_n)x_1^n)-A_n(x_1^n)$, it follows from \thetag{1-4} that
  $$
 E_n(I\times T_n\times I) E_n^{-1}= (\widehat{T}_n)_{\alpha_n}\times I_Z.
  $$ 
  Hence  $(\widehat{T}_n)_{\alpha_n}\times I_Z\in \Cal{G}$. Although, each $\alpha_n$ is a coboundary for the $\Bbb Z$-action given by $T_n$, the cocycle $\alpha$ is not a coboundary  for  $\Gamma$.
   In fact, we will now show the following.

{\bf Claim I.} {\it The cocycle $\alpha$  of $\Gamma $ is ergodic.}

For that we will show that each element $a\in G$ is an essential value of $\alpha$.
Given $n>0$ and a subset $B\subset X_1\times\cdots\times X_n$, denote by $[B]_1^n\subset X$ the corresponding cylinder with the ``head" $B$, i.e.
$[B]_1^n:=\{x\in X\mid x_1^n\in B\}$.
Let $U$ be a symmetric neighborhood  of $0$ in $G$.
Choose a countable partition $\Cal P$ of $G$ into Borel subsets $\Delta$ such that
$g-h\in U$ for all $g,h\in\Delta$ and each $\Delta\in\Cal P$.
Let $\psi_k:=\nu_k\circ f_k^{-1}$ for each $k>0$.
Using the ATI-assumption, we can find $m>n$ such that
$$
\|\psi_{n+1}*\cdots*\psi_m*\delta_a-\psi_{n+1}*\cdots*\psi_m\|_1<\epsilon.\tag1-5
$$
For each $\Delta\in\Cal P$, we let
$$
\aligned
A_\Delta&:=\Big\{y=(y_k)_{k=n+1}^m\in X_{n+1}\times\cdots \times X_m\,\Big|\,
\sum_{k=n+1}^mf_k(y_k)\in\Delta\Big\}\quad\text{and}\\
B_\Delta&:=\Big\{y=(y_k)_{k=n+1}^m\in X_{n+1}\times\cdots \times X_m\,\Big|\,
a+\sum_{k=n+1}^mf_k(y_k)\in\Delta\Big\}.
\endaligned
\tag1-6
$$
Then $\{A_\Delta\}_{\Delta\in\Cal P}$ and $\{B_\Delta\}_{\Delta\in\Cal P}$ are two measurable partitions of $X_{n+1}\times\cdots\times X_m$.
It follows from
\thetag{1-5} that 
$$
\align
\sum_{\Delta\in\Cal P}|\nu_{n+1}^m(A_\Delta)-\nu_{n+1}^m(B_\Delta)|&=
\sum_{\Delta\in\Cal P}|\psi_{n+1}*\cdots*\psi_m*\delta_a(\Delta)-\psi_{n+1}*\cdots*\psi_m(\Delta)|\\
&\le
\|\psi_{n+1}*\cdots\psi_m*\delta_a-\psi_{n+1}*\cdots *\psi_m\|_1\\
&<\epsilon,
\endalign
$$
where $\nu_{n+1}^m$ denotes the direct product $\bigotimes_{k=n+1}^m\nu_k$.
We can find subsets $A_\Delta'\subset A_\Delta$
and $B_\Delta'\subset B_\Delta$
such that
$$
\nu_{n+1}^m(A_\Delta')=
\nu_{n+1}^m(B_\Delta')=\min(\nu_{n+1}^m(A_\Delta),\nu_{n+1}^m(B_\Delta)).
\tag1-7
$$
Note that the group $\Gamma_{n+1,m}$ generated by $m-n$ mutually commuting transformations  $T_{n+1}\times I\times\cdots\times I$, $I\times T_{n+2}\times I\times\cdots\times I$, \dots, $I\times\cdots\times I\times T_m\in\text{Aut}_0(X_{n+1}\times\cdots\times X_m,\nu_{n+1}^m)$ is ergodic.
Hence, in view of \thetag{1-7},    Hopf's lemma \cite{HaOs}  yields that there is 
 a transformation $S_0\in[\Gamma_{n+1,m}]$ such that
 $S_0 A_\Delta'= B_{\Delta}'$
for each $\Delta\in\Cal P$.
We note that  
$$
\sum_{\Delta\in\Cal P}\nu_{n+1}^m(A_\Delta\setminus A_\Delta')
\le
\sum_{\Delta\in\Cal P}|\nu_{n+1}^m(A_\Delta)-\nu_{n+1}^m(B_\Delta)|<\epsilon.
$$
It follows  that  $\nu_{n+1}^m(\bigsqcup_{\Delta\in\Cal P}A_\Delta')>1-\epsilon.$
On the other hand, in view of \thetag{1-6},
for each $y\in A^+:= \bigsqcup_{\Delta\in\Cal P}A_\Delta'$,
$$
\bigg(\sum_{k=n+1}^mf_k\bigg)(y)-
\bigg(\sum_{k=n+1}^mf_k\bigg)(S_0 y)\in a+ U.
$$
We  now ``extend'' $S_0$ to a transformation $S\in\text{Aut}_0(X,\mu)$ by setting
 %$S\in\text{Aut}_0(X_1\times\cdots\cdots X_m,\bigotimes_{k=1}^m\nu_k)$ by setting
$$
Sx:=(x_1^n,S_0x_{n+1}^m, x_{m+1}^\infty)\in X\quad \text{for all }x\in X.
$$
Then $S\in[\Gamma]$ and 
$$
\alpha(S,x)\in a+ U\quad\text{whenever $x_{n+1}^m\in A^+$.}
\tag 1-8
$$
Then  we have that
$[B\times A^+]_1^m\subset [B]_1^n$,
$S[B\times A^+]_1^m\subset [B]_1^n$,
$\mu([B\times A^+]_1^m)>\frac 12\mu([B]_1^n)$ and \thetag{1-8} holds for all $x\in [B\times A^+]_1^m$.
Since the set of all cylinders is dense in $\goth B$, it follows from Lemma~C that $a$ is an essential value of $\alpha$.
Thus, Claim~I is proved.

To complete the proof of the theorem, we have already noticed that
 $(\widehat T_n)_{\alpha_n}\times I_Z\in\Cal G$ for each $n\in\Bbb N$.
 Hence $F(\gamma_\alpha(x,g),z)=F(x,g,z)$ at a.e. $(x,g,z)\in X\times G\times Z$ for each
 $\gamma\in\Gamma$.
 Claim~I yields that there is a function $M:Z\to\Bbb R$ such that $F(x,g,z)=M(z)$ at
 a.e. $(x,g,z)\in X\times G\times Z$.
 Since $F$ is invariant under $T_f\times C$, we obtain that $M$ is invariant under $C$.
 Since $C$ is ergodic, $M$ is constant a.e. and hence $F$ is constant a.e., i.e. $T_f\times C$ is ergodic.
 \qed
\enddemo

We call the cocycle $f$ in the statement of Theorem~1.5 {\it a cocycle of product type.}

\subhead 1.6. Application to Gaussian cocycles
\endsubhead
Let $(X,\goth B,\mu, T)$ be an ergodic Gaussian dynamical system.
It is completely determined by a restriction of the corresponding Koopman unitary operator $U_T$ to a closed (real) Gaussian subspace $H\subset L^2_0(X,\mu)$, called  {\it the first chaos}. (See, e.g. \cite{LePaTh} for the definitions.)
Let $\kappa$ denote the maximal spectral type of $U_T\restriction H$.
It is known that $T$ is ergodic if and only if $T$ is weakly mixing if and only if $\kappa$ is nonatomic.
Take $f\in H$.
Then the measurable map $f:X\to\Bbb R$
 considered as a cocycle of $T$ is called a {\it Gaussian} cocycle.
 It was shown in  \cite{LeLeSk} that if 
 $f$
 is a $T$-coboundary, i.e. $f=h\circ T-h$ for a measurable function $h:X\to\Bbb R$, then $h\in H$.
We now recall a conjecture from \cite{LeLeSk}.

\proclaim{Conjecture I} If a Gaussian cocycle $f$ is not a coboundary then $f$ is ergodic.
\endproclaim

We now prove this conjecture (in fact, we prove a stronger result) under an additional assumption that  $T$ is mildly mixing.

\proclaim{Theorem 1.6} If $T$ is a mildly mixing Gaussian transformation and $f$ is a Gaussian cocycle of $T$ which is not a coboundary then $T_f$
is sharply weak mixing.
\endproclaim

\demo{Proof}
Since $f\in H$, it follows  that  $\int_Xf\,d\mu=0$.
Hence, by  Atkinson's theorem \cite{At}, $T_f$ is conservative.
Consider  now the spectral decomposition for the pair $(H,U_T)$:
$$
H=\int^\oplus_\Bbb T\Cal H_z\,d\kappa(z)\quad\text{and}\quad
U_T=\int_{\Bbb T}^\oplus zI_z\,d\kappa(z),
$$
where $\Bbb T\ni z\mapsto\Cal H_z$ is the corresponding measurable field of Hilbert spaces and $I_z$ is the identity operator in $\Cal H_z$.
In other words, we can consider an element $h$ of $H$  as a measurable map $\Bbb T\ni z\mapsto h(z)\in\Cal H_z$ such that $\|h\|^2=\int_\Bbb T\|h(z)\|^2d\kappa(z)<\infty$.
We now let  $\Delta_n:=\big\{z\in\Bbb T\mid \frac 1{n+1}<|z-1|\le\frac 1{n}\big\}$.
Then we obtain a countable partition $\bigsqcup_{n=1}^\infty\Delta_n$
of $\Bbb T\setminus\{1\}$.
Since $\kappa(\{1\})=0$,
this countable partition generates a decomposition of $H$ into a direct sum $\bigoplus_{n\in\Bbb N}H_n$ of closed $U_T$-invariant subspaces $H_n$ consisting of the measurable maps $h:\Bbb T\ni z\mapsto h(z)\in \Cal H_z$ such that $h(z)=0$ whenever $z\not\in \Delta_n$.
This decomposition induces a decomposition of $(X,\mu, T)$ into the infinite direct product
$(X,\mu,T)=\bigotimes_{n=1}^\infty(X_n,\mu_n, T_n)$, where  $(X_n,\mu_n, T_n)$ is the Gaussian dynamical system associated with the pair $(H_n, U_T\restriction H_n)$ for each $n\in\Bbb N$.
Now we can expand $f$ into an orthogonal sum $f=\bigoplus_{n=1}^\infty f_n$ with $f_n\in H_n$ for each $n\in\Bbb N$.
Of course, for each $n>0$, there is $a_n\in H_n$ such that $f_n=U_Ta_n-a_n$.
Indeed, it follows from this equation that $f_n(z)=za_n(z)-a_n(z)$  and hence
$a_n(z)=(z-1)^{-1}f_n(z)$ for a.e. $z\in\Delta_n$.
Since $|z-1|^{-1}<n+1$ for all $z\in\Delta_n$, we obtain that  $a_n\in H_n$.
This yields an expansion 
$$
f=\bigoplus_{n=1}^\infty(U_Ta_n-a_n)=\bigoplus_{n=1}^\infty(a_n\circ T_n^{-1}-a_n)\tag1-9
$$ 
of $f$ into an infinite sum of $T_n$-coboundaries. 
Of course, $\sum_{n\in\Bbb N}\|a_n\|^2=+\infty$.
Otherwise the series $\sum_{n\in\Bbb N}a_n$ converges in $H$ and hence $f$ would be a coboundary  which contradicts  the assumption of the theorem.
We have that $\mu_n\circ a_n^{-1}=\Cal N_{0,\|a_n\|^2}$ for each $n\in\Bbb N$.
Passing, if necessary, to a subsequence we may assume without loss of generality that the convergence
in \thetag{1-9} is almost everywhere. 
Example~1.3 yields that the sequence $(\mu_n\circ a_n^{-1})_{n=1}^\infty$ is ATI.
It now follows from Theorem~1.5 that $T_f$ is sharply weak mixing. \qed
\enddemo

Consider now the general case.
Then there is a maximal (with respect to $\kappa$) subset $A$ of $\Bbb T$ such that
 $U_T$ restricted to the closed subspace $\int^\oplus_A\Cal H_zd\kappa(z)$ of $H$ is mildly mixing.
 We note that $A$ is symmetric.
 Then $\kappa$ decomposes into a sum of two orthogonal measures: $\kappa_{mm}:=\kappa\restriction A$ (the mildly mixing part of $\kappa$) and $\kappa_r:=\kappa\restriction (\Bbb T\setminus A)$ (the rigid part of $\kappa$).
 This decomposition defines a decomposition of $(X,\mu,T)$ into a direct product $(X_1,\mu_{mm},M)\times(X_2,\mu_{r},R)$, where $(X_1,\mu_{mm},M)$ is the Gaussian dynamical system corresponding to the pair $(\int^\oplus_A\Cal H_zd\kappa_{mm}(z), U_T)$
 and $(X_2,\mu_{r},R)$ is the Gaussian dynamical system corresponding to the pair $(\int^\oplus_{\Bbb T\setminus A}\Cal H_zd\kappa_r(z), U_T)$.
 Also, we obtain  a decomposition of $f$ into a sum $f_{mm}+f_r$, where $f_{mm}:=f1_A$ and $f_r=f1_{\Bbb T\setminus A}$.
 There are two possible cases: either $f_{mm}$ is a coboundary or 
 $f_{mm}$ is not a coboundary.
 In the first case $T_f$ is isomorphic to $Q\times R_{f_r}$.
 Moreover, $f_{r}$ is not a coboundary because otherwise $f$ would be  a coboundary.
 Since $Q$ is mildly mixing and $T_f$ is conservative, $T_f$ is ergodic if and only if  $R_{f_r}$ is ergodic.
 In the second case, $T_f$ is isomorphic to $Q_{f_{mm}}\times R_{f_r}$ and 
 $Q_{f_{mm}}$ is sharply weak mixing by Theorem~1.5.
 Since $T_f$ is conservative, it follows that $T_f$ is ergodic if and only if $R_{f_r}$ is ergodic.
 Thus, we have reduced the conjecture from \cite{LeLeSk} to the following one.
 
 \proclaim{Conjecture II} If a Gaussian cocycle $f$ is not a coboundary and $\kappa$ has only rigid part then
 $f$ is ergodic.\footnote{In \cite{LeLeSk}, there were constructed some concrete rigid Gaussian transformations admitting   ergodic Gaussian cocycles. In \cite{MaRa} this result was extended  to arbitrary rigid Gaussian transformations which have at least one Gaussian non-coboundary. However, it is unknown whether the ergodicity holds for each  Gaussian non-coboundary in those examples.}  
 \endproclaim

\head 2.  Krieger's type of  infinite direct products of dynamical systems of finite type \endhead

\subhead 2.1.  IDPFT systems
\endsubhead
Let $T_n$ be a  nonsingular invertible  transformation of a standard probability space $(X_n,\goth B_n,\mu_n)$ for each $n\in\Bbb N$.
Denote by $T$  the infinite direct product of $T_n$, $n\in\Bbb N$, acting
on the infinite product space $(X,\goth B,\mu):=\bigotimes_{n\in\Bbb Z}(X_n,\goth B_n,\mu_n)$.
By Theorem~A,
$T$ is $\mu$-nonsingular if and only if
$$
\prod_{n=1}^\infty \Big(1- H^2(\mu_n\circ T_n^{-1},\mu_n)\Big)>0\text{ \ or \  }\sum_{n=1}^\infty H^2(\mu_n\circ T_n^{-1},\mu_n)<\infty.\tag2-1
$$
If \thetag{2-1} does not hold then $\mu\circ T^{-1}\perp\mu$.
If $T$ is $\mu$-nonsingular then
$$
\frac{d\mu\circ T^{-1}}{d\mu}(x)=\prod_{n=1}^\infty
	\frac{d\mu_n\circ T_n^{-1}}{d\mu_n}(x_n)\qquad\text{at a.e. $x\in X$.}
	$$
Suppose now that $T_n$ is 
 {\it of finite type}, i.e that there exists a $\mu_n$-equivalent probability measure $\nu_n$ which is invariant under $T_n$ for each $n\in\Bbb N$.
We then put $\phi_n:=\frac{d\mu_n}{d\nu_n}$.
Since $1-H^2(\mu_n\circ T_n^{-1},\mu_n)=\int_{X_n}\sqrt{\frac{\phi_n\circ T_n^{-1}}{\phi_n}}\phi_nd\nu_n$,
the formula \thetag{2-1} and Theorem~A yield the following. 

\proclaim{Corollary  2.1}
$T$ is $\mu$-nonsingular if and only if
$$\prod_{n=1}^\infty\int_{X_n}
\sqrt{\phi_n\cdot \phi_n\circ T_n^{-1}}d\nu_n>0.%\tag1-7
\tag2-2
$$
$\mu\perp\nu$  if and only if 
$$
 \prod_{n=1}^\infty\int_{X_n}\sqrt \phi_nd\nu_n=0.
 %\tag1-8
 \tag2-3
$$
\endproclaim

 \definition{Definition 2.2} If $T$ is $\mu$-nonsingular and $T_n$ is of finite type for all $n>0$ then we say that the dynamical system $(X,\goth B,\mu,T)$ is {\it IDPFT (i.e. infinite direct product  of finite types)}.
 \enddefinition

Our purpose in this section is to investigate dynamical properties of IDPFT-systems.
The first result is about ergodicity of conservative IDPFT systems under the mild  mixing assumption on the factors.

\proclaim{Proposition 2.3}
Let $(X_n,\nu_n,T_n)$ be mildly mixing  for each $n>0$
and \thetag{2-2} and~\thetag{2-3} hold.
Suppose that $T$ is $\mu$-conservative.
Then  $T$ is $\mu$-sharply weak mixing and $\mu\perp\nu$.
\endproclaim

\demo{Proof}
Let $S$ be an ergodic conservative transformation of a standard probability space
$(Y,\goth C,\nu)$.
As in the proof of Theorem~1.5 one can show that $T\times S$ is either totally dissipative
or conservative.
Suppose  that $T\times S$ is conservative.
We have to prove that it is ergodic.
Let a subset $A\in\goth B\otimes\goth C$ be invariant under $T\times S$.
It follows from Theorem~B that for each $n>0$, 
$A$ belongs to the $\sigma$-algebra $\{\emptyset, X_1\times\cdots\times X_n\}\otimes\goth B_{n+1}\otimes\goth B_{n+2}\otimes\cdots\otimes\goth C$\footnote{When applying Theorem B, we consider the measure $(\bigotimes_{k=1}^n\nu_k)\otimes(\bigotimes_{k>n}\mu_k)\otimes\nu$ on $X\times Y$.
This measure is equivalent to $\mu\otimes\nu$.}.
By the Kolmogorov 0-1 law, the intersection of these $\sigma$-algebras is $\goth N\otimes\goth C$, where $\goth N$ is the trivial $\sigma$-algebra on $X$.
Thus $A=X\times D$ for some subset $D\in\goth C$.
Since $A$ in invariant under $T\times S$, it follows that $D$ is invariant under $S$.
Since $S$ is ergodic, we obtain that
 either $\mu\otimes\nu(A)=0$ or $\mu\otimes\nu(A)=1$.
\qed
\enddemo

\remark{Remark 2.4}
In Section~3 below
 we will give examples of $(X,\nu,T)$ and $\mu$ such that $(X,\mu,T)$ is of type $III_1$.
 In particular, there is no $\mu$-equivalent invariant probability measure.
 On the other hand, we do not know examples in which $(X,\mu,T)$ is of type $II_1$, i.e.
$T$ is mildly mixing with respect to a $\mu$-equivalent invariant probability measure.
  \endremark

\subhead 2.2. Radon-Nikodym cocycle and type $III_1$
\endsubhead
 Let $\Gamma$ be an ergodic countable subgroup of Aut$(Y,
\nu)$.
Denote by $\rho_\nu:\Gamma\times Y\to\Bbb R$ the logarithm of {\it the Radon-Nikodym}
cocycle of $\Gamma$, i.e.
$$
\rho_\nu(\gamma,y):=\log\frac{d\nu\circ \gamma}{d\nu}(y).
$$
 The $\rho_\nu$-skew product extension $\Gamma_{\rho_\nu}$ of $\Gamma$ is called the {\it Maharam extension} of $R$.
We note that $\Gamma_{\rho_\nu}$ preserves an equivalent $\sigma$-finite measure 
$\nu\otimes\kappa$, where 
 $\kappa$ is a Lebesgue absolutely continuous $\sigma$-finite measure on $\Bbb R$ such that $d\kappa(t)=e^{-t}dt$ for all $t\in\Bbb R$.
Similar to the finite measure preserving case, $\rho_\nu$ ``extends'' to the full group $[\Gamma]$
 in such a way that the cocycle identity holds.
 Moreover, we do not need the freeness condition for $\Gamma$ to define this extension. 
 
 We note that $\rho_\nu$ is a coboundary if and only if there is a $\Gamma$-invariant $\nu$-equivalent $\sigma$-finite measure on $(Y,\goth C)$.

 By the Maharam theorem (see \cite{Sc}), $\Gamma_{\rho_\nu}$ is conservative
 if and only if $\Gamma$ is conservative.
However  if $\Gamma$ is ergodic then 
  $\Gamma_{\rho_\nu}$ is not necessarily ergodic.
If the Maharam extension of $\Gamma$ is ergodic then
$\Gamma$ is called {\it of Krieger's type $III_1$}.
If for each homomorphism $\vartheta:\Gamma\to\text{Aut}_0(Y,\nu)$ such that the image $\{\vartheta(\gamma)\mid\gamma\in\Gamma\}$ is ergodic, the direct product
$\{\gamma\times\vartheta(\gamma)\mid\gamma\in\Gamma\}$ is ergodic and of type $III_1$ then $\Gamma$ is said {\it to be of  stable Krieger's type $III_1$}.

It is possible to define essential values of $\rho_\nu$ in the same way as in the finite measure preserving case.

\definition{Definition 2.5} 
An element $g\in \Bbb R$ is called an {\it essential value of $\rho_\nu$}
if for each subset $A\subset Y$ of positive measure and a neighborhood $U$ of $g$, there are a Borel subset $B\subset A$ and an element $\gamma\in\Gamma$ such that $\nu(B)>0$,
$\gamma B\subset A$ and $\rho_\nu(\gamma, y)\in U$ for all $y\in B$.
\enddefinition

We refer to \cite{Sc} and \cite{HaOs} for the proof of the following results:
\roster
\item"---"
The set $r(\rho_\nu)$ is a closed subgroup in $\Bbb R$.
\item"---"
$\Gamma$ is of type $III_1$ if and only if $r(\rho_\nu)=\Bbb R$.
\endroster

If there is $\lambda\in(0,1)$ such that $r(\rho_\nu)=\{n\log\lambda\mid n\in\Bbb Z\}$ then
 $\Gamma$ is said to be {\it of Krieger's  type $III_\lambda$.}

We will need the following analog of Lemma~C.
It follows from a more general \cite{Da, Lemma 2.4}.
 
 \proclaim{Lemma D}
 An element $a\in\Bbb R$ is an essential value of $\rho_\nu$ if there exists $\delta>0$ such that for each $\epsilon>0$ and each subset $B$ from a dense collection $\goth C_0$ of subsets in $\goth C$, there is a subset $B_0\subset B$ and a transformation
 $\theta\in[\Cal \Gamma]$
 such that $\nu(B_0)>\delta\nu(B)$, $\theta B_0\subset B$ and either $|\rho_\nu(\theta, y)-a|\le\epsilon$ for all $y\in B_0$ or $|\rho_\nu(\theta, y)+a|\le\epsilon$ for all $y\in B_0$.
 \endproclaim

 \subhead 2.3. On conservativeness of IDPFT systems
 \endsubhead
 In this subsection we first establish a general result on conservativeness of infinite direct product systems.
 
 \proclaim{Proposition 2.6} Let $(X_n,\goth B_n,\mu_n, T_n)$ be an ergodic nonsingular dynamical system on a standard probability space for each $n\in\Bbb N$ and let \thetag{2-1} hold.
 Let $(X,\goth B,\mu, T):=\bigotimes_{n=1}^\infty(X_n,\goth B_n,\mu_n, T_n)$.
 If, for each $n\in\Bbb N$, there is a function $\alpha_n:X_n\to[1,+\infty)$ such that for each $k\in\Bbb N$
 $$
\alpha_n(x)^{-1}\le\frac{d\mu_n\circ T^k}{d\mu_n}(x)\le\alpha_n(x)\quad\text{ at a.e. $\mu_n$-a.e. $x\in X_n$}
 $$
  then the dynamical system $(X,\goth B,\mu, T)$ is either conservative or totally dissipative.
 Moreover, if $(Y,\goth C,\nu, S)$ is an ergodic conservative nonsingular dynamical system then the direct product
 $T\times S$ is either conservative or totally dissipative.
 \endproclaim
 
 \demo{Proof} We will prove the second claim only.
 By the Hopf criterion \cite{DaSi, \S 2}, 
 $$
 \Cal D(T\times S)=\bigg\{(x,y)\in X\times Y\mid \sum_{k=1}^\infty\frac{d(\mu\otimes\nu)\circ (T\times S)^k}{d(\mu\otimes\nu)}(x,y)<\infty\bigg\}.
 $$
 For each $r>0$, we consider a transformation $\gamma_r$ of $X$ by setting
 $\gamma_r(x_1,x_2,\dots):=(x_1,\dots,x_{r-1},T_rx_r,x_{r+1},\dots)$.
 Of course, $\gamma_r\in\text{Aut}(X,\mu)$.
 Denote by $\Gamma$ the transformation group generated by $\gamma_r$, $r\in\Bbb N$.
It follows from the Kolmogorov 0-1 law that $\Gamma$ is ergodic. 
 We claim that $ \Cal D(T\times S)$ is invariant under $\gamma_r\times I$
 for each $r$.
Let $(x,y)\in\Cal D(T\times S)$.
Since for each $k>0$,
$$
 \align
 \frac{d\mu\circ T^k}{d\mu}(\gamma_rx)&= \frac{d\mu_r\circ T_r^k}{d\mu_r}(T_rx_r)
 \bigg( \frac{d\mu_r\circ T_r^k}{d\mu_r}(x_r)\bigg)^{-1}
  \prod_{n=1}^\infty \frac{d\mu_n\circ T_n^k}{d\mu_n}(x_n)\\
  &\le
  \alpha_r(T_rx_r)\alpha_r(x_r)^{-1}
   \frac{d\mu\circ T^k}{d\mu}(x),
  \endalign
 $$
 it follows that
 $$
 \align
 \sum_{k=1}^\infty\frac{d(\mu\otimes\nu)\circ (T\times S)^k}{d(\mu\otimes\nu)}(\gamma_rx,y)
 &=
 \sum_{k=1}^\infty \frac{d\mu\circ T^k}{d\mu}(\gamma_rx)\frac{d\nu\circ S^k}{d\nu}(y)\\
& \le
\frac{\alpha_r(T_rx_r)}{\alpha_r(x_r)}\sum_{k=1}^\infty \frac{d\mu\circ T^k}{d\mu}(x)\frac{d\nu\circ S^k}{d\nu}(y)\\
&=\frac{\alpha_r(T_rx_r)}{\alpha_r(x_r)}
\sum_{k=1}^\infty\frac{d(\mu\otimes\nu)\circ (T\times S)^k}{d(\mu\otimes\nu)}(x,y)<\infty.
\endalign
  $$
 Thus, $(\gamma_rx,y)\in\Cal D(T\times S)$.
 Since $\Cal D(T\times S)$ is invariant under $I\times S$, we obtain that  $\Cal D(T\times S)$ is invariant under an ergodic transformation group on $X\times Y$ generated by $I\times S$ and $\gamma\times I$, $\gamma\in\Gamma$.
 Hence,  either $(\mu\otimes\nu)(\Cal D(T\times S))=0$ or 
 $(\mu\otimes\nu)(\Cal D(T\times S))=1$, as desired. \qed
  \enddemo
 
 We now apply Proposition 2.6 to IDPFT systems.
 
 \proclaim{Corollary 2.7} Let $(X_n,\goth B_n,\mu_n, T_n)$ be an ergodic nonsingular dynamical system on a standard probability space for each $n\in\Bbb N$ and let \thetag{2-1} hold.
 Suppose that for each $n\in\Bbb N$, there is a $\mu_n$-equivalent $T_n$-invariant
 probability measure on $X_n$.
 Let $(X,\goth B,\mu, T):=\bigotimes_{n=1}^\infty(X_n,\goth B_n,\mu_n, T_n)$.
 Then the dynamical system $(X,\goth B,\mu, T)$ is either conservative or totally dissipative.
 Moreover, if $(Y,\goth C,\nu, S)$ is an ergodic conservative nonsingular dynamical system then the direct product
 $T\times S$ is either conservative or totally dissipative.
  \endproclaim
 \demo{Proof}
 Let $\phi_n:=\frac{d\mu_n}{d\nu_n}$ for each $n>0$.
 If for each $n>0$, there is a real $\alpha_n\ge 1$ such that $\alpha_n^{-1}\le\phi_n\le\alpha_n$ almost everywhere then the claim of the corollary follows directly from Proposition~2.6.
 We now show that the general case can be reduced to the ``bounded'' one.
 Indeed, for each $n>0$, we can find a probability measure $\widetilde\mu_n\sim\mu_n$
 such that $H^2(\widetilde\mu_n,\mu_n)\le 2^{-n}$ and the Radon-Nikodym derivative
 $\frac{d\widetilde\mu_n}{d\nu_n}$ is bounded from above and separated from $0$ from below.\footnote{For that, take the Radon-Nikodym derivative $\frac{d\mu_n}{d\nu_n}$ and change it on a subset of very small measure to get the boundedness. The ``modified'' function will be the Radon-Nikodym derivative $\frac{d\widetilde\mu_n}{d\mu_n}$.}
 Since
 $\sum_{n=1}^\infty H^2(\widetilde\mu_n,\mu_n)<\infty$, it follows from
Theorem~A that $\mu\sim\widetilde\mu:=\bigotimes_{n=1}^\infty\widetilde\mu_n$.
It remains to note that the conservativeness of a dynamical system does not depend on the choice of quasi-invariant measure within its equivalence class.
\qed
  \enddemo

 \subhead 2.4. Sharp weak mixing for Maharam extensions of IDPFT systems
 \endsubhead
 We first  introduce a ``nonsingular analog'' of the property ATI (cf. Definition~1.2).
 
 \definition{Definition 2.8}
A sequence $(\xi_n)_{n=1}^\infty$ of probability non-atomic Borel measures on $G$ is called
asymptotically translation quasi-invariant (ATQI) if for each $a\in G$ there exists $\zeta_a>0$ such that
 for every $n\in\Bbb N$  there are $m>n$   and a Borel subset $W_{n,m}\subset G$ such that
$$
\gather
 \zeta_a\le (\xi_n*\xi_{n+1}*\cdots*\xi_{n+m})(W_{n,m}),\\
\xi_n*\xi_{n+1}*\cdots*\xi_{n+m}*\delta_a\prec\xi_n*\xi_{n+1}*\cdots*\xi_{n+m}\quad\text{and}\\
\zeta_a\le\frac{d(\xi_n*\xi_{n+1}*\cdots*\xi_{n+m}*\delta_a)}{d(\xi_n*\xi_{n+1}*\cdots*\xi_{n+m})}(t)\quad\text{for each $t\in W_{n,m}$}.
\endgather
$$ 
\enddefinition

We will need the following lemma on continuous measures.
 
 \proclaim{Lemma 2.9} Given a standard  probability space $(Y,\goth C,\nu)$, a nonnegative function $\phi\in L^1(Y,\nu)$ such that the measure $\nu\circ \phi^{-1}$ is nonatomic
 and $\delta\in (0,1)$, then
 $$
 \max\bigg\{\int_A\phi \,d\nu\,\bigg|\, \nu(A)=\delta\bigg\}\ge\frac{\delta}2\int_Y\phi\,d\nu.
 $$
 \endproclaim
 \demo{Proof}
Since $\nu\circ \phi^{-1}$ is nonatomic, 
the function
 $$
 f:[0,+\infty)\ni t\mapsto f(t):=\nu(\{y\in Y\mid \phi(y)>t\})\in [0,1]
 $$
 is continuous and non-increasing.
 Since $f(0)=1$ and $\lim_{t\to\infty}f(t)=0$,
 there is   a unique $t_0>0$ such that $f(t_0)=\delta$.
 Let $E:=\{y\in Y\mid \phi(y)>t_0\}$.
 Then $\max\{\int_A\phi d\nu\mid \nu(A)=\delta\}=\int_E\phi\,d\nu$.
 Find $n\ge 1$ such that $\frac 1{n+1}<\delta\le \frac1n$.
Then there is a partition  $Y=Y_1\sqcup\dots\sqcup Y_{n+1}$  of $Y$ into subsets $Y_k$ such that
$\nu(Y_k)=\delta$ for each $k=1,\dots,n$ and $\mu(Y_{n+1})\le\delta$.
We now have:
$$
\int_Y\phi\,d\nu= \sum_{k=1}^{n+1}\int_{Y_k}\phi\,d\nu\le (n+1)\int_E\phi\,d\nu\le \frac {n+1}{n\delta}
\int_E\phi\,d\nu\le \frac 2{\delta}\int_E\phi\,d\nu.
$$
 \qed
 
 \enddemo

The next theorem is a nonsingular analogue of Theorem~1.5.
The skeleton of the proof is similar to that of Theorem~1.5.

 \proclaim{Theorem 2.10}
 Let a dynamical system $(X_n,\goth B_n,\nu_n,T_n)$ be mildly mixing  for each $n>0$.
 Let $\mu_n$ be a probability on $X_n$ such that $\mu_n\sim\nu_n$ for each $n\in\Bbb N$.
Let $\phi_n:=\frac{d\mu_n}{d\nu_n}$ and  \thetag{2-2} hold. 
We set
$$
(X,\goth B,\nu, T):=\bigotimes_{n\in\Bbb N}(X_n,\goth B_n,\nu_n, T_n)
$$
 and $\mu:=\bigotimes_{n=1}^\infty\mu_n$.
 If $T$ is $\mu$-conservative and the sequence of probability measures $(\nu_n\circ(\log\phi_n)^{-1})_{n=1}^\infty$ is ATQI then $T\in\text{\rom{Aut}}(X,\mu)$ is ergodic of stable type $III_1$.
 Moreover, the Maharam extension of $T$ is sharply weak mixing.
 \endproclaim

\demo{Proof} By the  Maharam theorem, the Maharam extension $T_{\rho_\mu}$ is conservative.
Let $C$ be  an ergodic conservative transformation of a standard probability space $(Z,\goth Z,\eta)$.
As in the proof of Theorem~1.5, one can show that $T_{\rho_\mu}\times C$
is either totally dissipative or conservative.
Suppose that it is conservative and prove that it is ergodic.

  Let a function $F\in L^\infty(X\times \Bbb R\times Z,\mu\otimes\kappa\otimes\eta)$ be invariant under $T_{\rho_\mu}\times C$.
  We first show that $F$ is also invariant under a huge group of transformations.
Fix $n>0$.
%For each $x\in X$, we write $x_1^n:=(x_1,\dots, x_n)\in X_1\times\cdots \times X_n$ and $x_{n+1}^\infty:=(x_{n+1},x_{n+2}, \dots)\in X_{n+1}\times X_{n+2}\times\cdots$.
%Then $x=(x_1^n,x_{n+1}^\infty)$.
We define a nonsingular automorphism $T^{(n)}$ of $\big(\bigotimes_{k>n}X_k,\mu^{(n)}\big)$, where $\mu^{(n)}:=\bigotimes_{k=n+1}^\infty\mu_k$,
and
a    measure preserving isomorphism  $E_n$  of $(X\times \Bbb R\times Z,\mu\otimes\kappa\otimes\eta)$ onto the product space 
$\big(X\times \Bbb R\times Z,\big(\bigotimes_{k=1}^n\nu_k\big)\otimes\mu^{(n)}\otimes\kappa\otimes\eta\big)$
by setting
$$
\align
T^{(n)}(x_k)_{k=n+1}^\infty&:=(T_{k}x_k)_{k=n+1}^\infty\quad\text{and}\\
E_n(x,t,z)&:=\Big(x,t+\sum_{k=1}^n\log\phi_k(x_k),z\Big).
%V_n((x_k)_{k=n+1}^\infty,g,z)&:=\Big((T_{k}x_k)_{k=n+1}^\infty, g+\sum_{k>n}(f_k(T_kx_k)-f_k(x_k)),Cz\Big).
\endalign
$$
Since
$$
\align
T_{\rho_\mu}(x,t)&=\bigg(T_1x_1,T_2x_2,\dots,t+\sum_{k=1}^\infty\log\frac{d\mu_k\circ T_k}{d\mu_k}(x_k)\bigg)\\
&=\bigg(T_1x_1,T_2x_2,\dots,t+\sum_{k=1}^\infty\Big(\log\phi_k( T_kx_k)-\log\phi_k(x_k)\Big)\bigg),
\endalign
$$
it follows that
$$
E_n (T_{\rho_\mu}\times C)E_n^{-1}=(T_1\times\cdots \times T_n)\times (T^{(n)})_{\rho_{\mu^{(n)}}}\times C.\tag2-4
$$
Since $T_{\rho_\mu}\times C$ is conservative, it follows from \thetag{2-4}
that the product $(T^{(n)})_{\rho_{\mu^{(n)}}}\times C$ is also conservative.
On the other hand,  the function $F\circ E_n^{-1}$ is invariant under 
$E_n(T_f\times C)E_n^{-1}$.
Utilizing these two facts plus the mild mixing of the transformation  $T_1\times\cdots \times T_n$,
we deduce  from Theorem~B
that $F\circ E_n^{-1}$ does not depend on the coordinates $x_1,\dots,x_n$.
Hence, for each transformation $S\in\text{Aut}_0(X_1\times\cdots\times X_n,\bigotimes_{k=1}^n\nu_k)$,
we have that
 $
 F\circ E_n^{-1}\circ (S\times I\times I_Z)=F\circ E_n^{-1}.
 $
 Therefore
$F$ is invariant under the transformation $E_n^{-1}(S\times I\times I_Z) E_n\in\text{Aut}_0(X\times \Bbb R\times Z,\mu\otimes\kappa\otimes\eta)$ and
$$
E_n^{-1}(S\times I\times I_Z) E_n=(S\times I)_{\rho_\mu}\times I_Z.\tag 2-5
%(S x_1^n,x_{n+1}^\infty, g-A_n(x^n_1)+A_n(S x_{1}^n), z),\tag1-4
$$
Denote by $\Gamma$ the group of  nonsingular transformations of   $(X,\goth B,\mu)$
 generated by $I\times T_n\times I$, $n\in\Bbb N$.
Then $\Gamma$ is  an ergodic Abelian countable  subgroup of Aut$(X,\mu)$ and
$F$ is invariant under $\{\gamma_{\rho_\mu}\times I_Z\mid \gamma\in\Gamma\}$ by \thetag{2-5}.
\comment

Thus, we have shown that $F$ is invariant under each transformation from the  set
$$
\Cal G:=\bigcup_{n>0}E_n\bigg(\text{Aut}_0\bigg(X_1\times\cdots\times X_n,\bigotimes_{k=1}^n\nu_k\bigg)\times \{I\}\bigg) E_n^{-1}.%\subset\text{Aut}_0(X,\nu).
$$

We now consider a new dynamical system.

Let a function $F\in L^\infty(X\times\Bbb R,\mu\times\kappa)$ be invariant under $\widetilde T$.
Fix $n>0$.
We define a  Borel map
$A_n:X_1\times\cdots\times X_n\to\Bbb R$ 
by setting 
$$
A_n(x_1,\dots,x_n):=\log\phi_1(x_1)+\cdots+\log\phi_n(x_n).
$$
Then that
the mapping
$C_n:X\times \Bbb R\to X\times R$,
given by
$$
C_n(x,t):=(x,t+A_nx_1^n),
$$
 is a measure preserving isomorphism
of $(X\times\Bbb R,\mu\otimes\kappa)$ onto the product space $(X\times\Bbb R,(\bigotimes_{k=1}^n\nu_k\otimes\bigotimes_{k=n+1}^\infty\mu_k)\otimes\kappa)$.
A direct verification shows that
$$
C_n\widetilde TC_n^{-1}=(T_1\times\cdots T_n)\times \widetilde{\bigotimes_{k=n+1}^\infty T_k}.
$$
We note that $F\circ C_n^{-1}$ is invariant under 
$C_n\widetilde TC_n^{-1}$.
The transformation $\bigotimes_{k=n+1}^\infty T_k$ is conservative as a factor of $T$.
By the Maharam theorem,  the Maharam extension of $\bigotimes_{k=n+1}^\infty T_k$ is also conservative.
By~Theorem~B,
$F\circ C_n^{-1}$ does not depend on $x_1,\dots,x_n$.

Hence, for each transformation $\tau\in\text{Aut}_0(X_1\times\cdots\times X_n,\bigotimes_{k=1}^n\nu_k)$,
we have that
 $F\circ C_n^{-1}\circ (\tau\times I)=F\circ C_n^{-1}$.
 Therefore
$F$ is invariant under the transformation $C_n^{-1}(\tau\times I) C_n\in\text{Aut}_0(X\times\Bbb R,\mu\times\kappa)$, and
$$
C_n^{-1}(\tau\times I) C_n(x,t)=
(\tau x_1^n,x_{n+1}^\infty, t+A_n(x^n_1)-A_n(\tau x_{1}^n))=
(\tau x_1^n,x_{n+1}^\infty, t+\frac{d\mu\circ\tau}{d\mu}(x)).
$$

In other words,  $C_n^{-1}(\tau\times I) C_n=\widetilde{\tau\times I}$.
In particular, we obtain that $F$ is invariant under the Maharam extension of every transformation from the group
$$
\Cal G:=\bigcup_{n>0}\bigg\{\tau\times I\,\bigg |\, \tau\in
%\text{Aut}_0\bigg(X_1\times\cdots\times X_n,\bigotimes_{k=1}^n\nu_k\bigg)
[T_1,\dots,T_n]
\bigg\}\subset\text{Aut}(X,\mu).
$$

\endcomment

{\bf Claim II}.  We claim that $\Gamma$ is of type $III_1$.
Equivalently,
we will show that
 each $a\in\Bbb R$ is an essential value for the cocycle $\rho_\mu$ of $\Gamma$.
For that, fix $n>0$, $\epsilon>0$
and a Borel subset $B\subset X_1\times\cdots\times X_n$.
Denote by $\psi_k$ the pushforward of $\nu_k$ under $\log\phi_k$ for each $k>0$.
By ATQI, there is $\zeta_a>0$ (which does not depend on $n$), 
 $m>n$ and a subset  $W_{n+1,m}\subset \Bbb R$
 such that
$$
\gathered
\zeta_a\le (\psi_{n+1}*\cdots*\psi_m)(W_{n+1,m}),\\
\psi_{n+1}*\cdots*\psi_m*\delta_a\prec\psi_{n+1}*\cdots*\psi_m\quad\text{and}\\
\zeta_a\le\frac{d(\psi_{n+1}*\cdots*\psi_{n+m}*\delta_a)}{d(\psi_{n+1}*\cdots*\psi_{n+m})}(t)\quad\text{for each $t\in W_{n+1,m}$}.
\endgathered
\tag2-6
$$
Choose a countable partition $\Cal P$ of $W_{n+1,m}$ into subsets of diameter no more than $\epsilon$.
For each $\Delta\in\Cal P$, we let
$$
\align
A_\Delta&:=\bigg\{y=(y_k)_{k=n+1}^m\in X_{n+1}\times\cdots \times X_m
\,\bigg|\,
\sum_{k=n+1}^m\log\phi_k(y_k)\in\Delta\bigg\}\qquad\text{and}\\
B_\Delta&:=\bigg\{y=(y_k)_{k=n+1}^m\in X_{n+1}\times\cdots \times X_m
\,\bigg|\,
a+\sum_{k=n+1}^m \log\phi_k(y_k)\in\Delta
\bigg\}.
\endalign
$$
%Then $\{A_\Delta\}_{\Delta\in\Cal P}$ and $\{B_\Delta\}_{\Delta\in\Cal P}$ are two measurable partitions of $W_{n+1,m}$.
%$X_{n+1}\times\cdots\times X_m$.
Let $\mu_{n+1}^m:=\bigotimes_{k=n+1}^m\mu_k$, $\nu_{n+1}^m:=\bigotimes_{k=n+1}^m\nu_k$ and $\phi^m_{n+1}:=\frac{d\mu_{n+1}^m}{d\nu_{n+1}^m}$.
Dropping off some atoms of $\Cal P$ if necessary,  we may assume without loss of generality  that $\nu_{n+1}^m(A_\Delta)>0$ (and hence
$\nu_{n+1}^m(B_\Delta)>0$ in view of \thetag{2-6}) for each $\Delta\in\Cal P$.
Note that the group $\Gamma_{n+1,m}$ generated by $m-n$ mutually commuting transformations  $T_{n+1}\times I\times\cdots\times I$, $I\times T_{n+2}\times I\times\cdots\times I$, \dots, $I\times\cdots\times I\times T_m\in\text{Aut}_0(X_{n+1}\times\cdots\times X_m,\nu_{n+1}^m)$ is ergodic.
Suppose that  $\nu_{n+1}^m(A_\Delta)>\nu_{n+1}^m(B_\Delta)$ for some $\Delta\in\Cal P$.
We note that for each  Borel subset $A'\subset A_\Delta$,
$$
\mu_{n+1}^m(A')
=\int_{A'}\nu_{n+1}^m(A_\Delta)\phi^m_{n+1}d\bigg(\frac{\nu^m_{n+1}}{\nu^m_{n+1}(A_\Delta)}\bigg).\tag2-7
$$
We now apply Lemma 2.9 to the space $A_\Delta$ equipped with the conditional measure $\nu^m_{n+1}(\cdot)/\nu^m_{n+1}(A_\Delta)$, the functionÊ $\nu_{n+1}^m(A_\Delta)\phi^m_{n+1}$ and $\delta=\nu^m_{n+1}(B_\Delta)/\nu^m_{n+1}(A_\Delta)$ to obtain a Borel subset 
 $A'_\Delta$ of conditional measure $\delta$, so
 $$
 \nu^m_{n+1}(A'_\Delta)=\nu^m_{n+1}(B_\Delta),
 $$
 and moreover (using~\thetag{2-7})
$$
\mu_{n+1}^m(A_\Delta')\ge \frac 12 \frac{\nu_{n+1}^m(B_\Delta)}{\nu_{n+1}^m(A_\Delta)}\mu_{n+1}^m(A_\Delta)\ge\frac{\zeta_a}2\,\mu_{n+1}^m(A_\Delta).
$$
%$$
%\mu_{n+1}^m(A_\Delta')\ge  \frac{\nu_{n+1}^m(A_\Delta')}{\nu_{n+1}^m(A_\Delta)}\mu_{n+1}^m(A_\Delta)\ge\zeta_a\,\mu_{n+1}^m(A_\Delta).
%$$
The latter inequality follows from \thetag{2-6} because 
$$ 
\frac{\nu_{n+1}^m(B_\Delta)}{\nu_{n+1}^m(A_\Delta)}=\frac{(\psi_n*\psi_{n+1}*\cdots*\psi_{n+m}*\delta_a)(\Delta)}{(\psi_n*\psi_{n+1}*\cdots*\psi_{n+m})(\Delta)}.
$$
By Hopf's lemma, there is 
 a transformation $S_0\in[\Gamma_{n+1,m}]$ such that
 \roster
 \item"---"
$S_0A_\Delta\subset B_\Delta$ if $\nu_{n+1}^m(A_\Delta)\le\nu_{n+1}^m(B_\Delta)$
and 
\item"---" $S_0A_\Delta'\subset B_\Delta$ if $\nu_{n+1}^m(A_\Delta)>\nu_{n+1}^m(B_\Delta)$.
\endroster
Let 
$$
A^+:=\bigsqcup_{\nu_{n+1}^m(A_\Delta)\le\nu_{n+1}^m(B_\Delta)}A_\Delta\sqcup\bigsqcup_{\nu_{n+1}^m(A_\Delta)>\nu_{n+1}^m(B_\Delta)}A_\Delta'.
$$
Then $\mu_{n+1}^m(A^+)\ge\frac{\zeta_a}2\mu_{n+1}^m(\bigsqcup_{\Delta\in\Cal P}A_\Delta)\ge \frac{\zeta_a^2}2$.
\comment
Let 
$$
A^-:=\bigsqcup_{\nu_{n+1}^m(A_\Delta)\le\nu_{n+1}^m(B_\Delta)}A_\Delta\quad\text{and}\quad A^+:=\bigsqcup_{\nu_{n+1}^m(A_\Delta)>\nu_{n+1}^m(B_\Delta)}A_\Delta.
$$
Then 
$$
\mu_{n+1}^m(A^-\sqcup A^+)=\nu_{n+1}^m(A^-\sqcup A^+)=1.
$$
We will  consider only the case where $\mu_{n+1}^m(A^-)\ge 1/2$.
(The other case, where $\mu_{n+1}^m(A^+)\ge1/2$ is considered in a similar way by changing the roles of $A_\Delta$ and $B_\Delta$.)

Since $\nu_{n+1}^m(B_\Delta)\ge \nu_{n+1}^m(A_{\Delta})$ for each $A_\Delta\subset A^-$,
   Hopf's lemma \cite{HaOs}  yields that there is 
 a transformation $S_0\in[\Gamma_{n+1,m}]$ such that
 $S_0 A_\Delta\subset B_{\Delta}$
whenever $A_\Delta\subset A^+$.

Then we can find 
a $\nu_{n+1}^m$-preserving invertible
transformation $S_0\in[T_{n+1},\dots, T_m]$ of $X_{n+1}\times\cdots\times X_m$
 such that

and  \roster
 \item
 if $\nu(A_\Delta)>\nu(B_\Delta)$
 then $\mu(A_\Delta')\ge\mu(A_\Delta\setminus A_\Delta')$,
 \item
$(\bigotimes_{k=n+1}^m\nu_k)(A_\Delta')=
(\bigotimes_{k=n+1}^m\nu_k)(B_\Delta')=\min(\bigotimes_{k=n+1}^m\nu_k)(A_\Delta),\bigotimes_{k=n+1}^m\nu_k)(B_\Delta))$
and
\item
 $S_0 A_\Delta'= B_{\Delta}'$
for each $\Delta\in\Cal P$.
\endroster
$S_0 A_\Delta\subset B_\Delta$ whenever $A_\Delta\subset A^-$.

\endcomment
Of course, for each $y\in A^+$,
$$
\bigg(\sum_{k=n+1}^m\log\phi_k\bigg)(y)-
\bigg(\sum_{k=n+1}^m\log\phi_k\bigg)(S_0 y)=a\pm\epsilon.\tag2-8
$$
We  now ``extend'' $S_0$ to a transformation $S\in\text{Aut}(X,\mu)$ by setting
 %$S\in\text{Aut}_0(X_1\times\cdots\cdots X_m,\bigotimes_{k=1}^m\nu_k)$ by setting
$$
S:=I\times S_0\times I.
$$
Then $S\in[\Gamma]$ and in view of \thetag{2-8}, 
$$
\rho_\mu(S,x)= -a\pm\epsilon\quad\text{whenever $x_{n+1}^m\in A^+$.}
\tag 2-9
$$
We now have that
$$
[B\times A^+]_1^m\subset [B]_1^n,\quad
S[B\times A^+]_1^m\subset [B]_1^n,\quad
\mu([B\times A^+]_1^m)\ge\frac{ \zeta_a^2}2\mu([B]_1^n)
$$
 and \thetag{2-9} holds for all $x\in [B\times A^+]_1^m$.
Since the set of all cylinders is dense in $\goth B$, it follows from Lemma~D that $a$ is an essential value of $\alpha$.
Thus, Claim~II is proved.

The assertion of the theorem follows from Claim II in the very same way as the assertion of 
Theorem~1.5 follows from Claim I (in the proof of Theorem~1.5).
\qed
\enddemo

\remark{Remark 2.11} In this remark we clarify some subtle points in the proof of Theorem~2.10.
Let $\goth B_0$ stand for the collection of all cylinders in $X$.
Then $\goth B_0$ is dense in $\goth B$ as with respect to $\mu$ as with respect to $\nu$.
Though $\mu\perp\nu$, the two measures are equivalent on $\goth B_0$, i.e.
$\mu(B)=0$ if and only if $\nu(B)=0$ whenever $B\in\goth B_0$.
Given a transformation $\theta\in\text{Aut}(X,\mu)$,
 the Radon-Nikodym derivative $\frac{d\mu\circ\theta}{d\mu}$ is defined up to a subset of zero $\mu$-measure.
Hence it has no sense as a function on $(X,\nu)$.
However, if we consider transformations of a specific product structure, say  $\gamma\in \Gamma$, then $\frac{d\mu\circ\gamma}{d\mu}$ is defined, in fact,  up
a subset of zero $\mu$-measure from $\goth B_0$.
Therefore, $\frac{d\mu\circ\gamma}{d\mu}$  is  well defined as a measurable function on $(X,\nu)$ as well.
Thus, the cocycle $\rho_\mu:\Gamma\times X\to\Bbb R$ is well defined simultaneously on $(X,\mu)$
and on $(X,\nu)$.
Another observation is that given a transformation $S_0\in[\Gamma_{n+1,m}]$, the extension
$S:=I\times S_0\times I$ of $S_0$ to $X$ is a well defined transformation  from Aut$(X,\mu)$ as well as from Aut$_0(X,\nu)$.
Thus, though an element of the full group $[\Gamma]$ is defined up to subset of zero measure, ``belonging'' $S\in[\Gamma]$ is well defined with respect to $\mu$ as well as with respect to $\nu$.
\endremark

\subhead 2.5. On type $III_\lambda$ for rigid IDPFT systems
\endsubhead
We would like to emphasize that
the conclusion of Theorem~2.10 does not hold if we drop the mild mixing condition on $T_n$ and
 the ATQI property.
We illustrate this on a family of IDPFT systems consisting of infinite product of periodic transformations.
Let $P=\{p_n\mid n\in\Bbb N\}$ be a countable subset of mutually coprime positive integers
and $p_1<p_2<\cdots$.
Of course, $p_n\to+\infty$ as $n\to\infty$.
In the examples that we are going to construct in this subsection, $p_n\to +\infty$ very fast (to be specified below).
For $n\in\Bbb N$, we set $X_n:=\{0,1,\dots,p_n-1\}$ and identify $X_n$ with the cyclic group $\Bbb Z/p_n\Bbb Z$.
Then $T_n:X_n\to X_n$, given by $T_nx=x+1\pmod {p_n}$, is a bijection of $X_n$.
The infinite product $T=\bigotimes_{n=1}^\infty T_n$ is a minimal rotation on the compact totally disconnected Abelian group $X:=\bigotimes_{n=1}^\infty X_n$.
Of course, the Haar measure $\nu$ on $X$ is the only $T$-invariant Borel probability measure on
$X$.
This measure is the infinite direct product of the equidistributions on $X_n$, $n\in\Bbb N$.
Fix $\lambda\in(0,1)$.
Let
$$
Y_n:=\bigg\{x_n\in X_n\mid x_n<\frac{p_n}2-\prod_{k=1}^{n-1}p_k\bigg\}\cup
\bigg\{x_n\in X_n\mid \frac{p_n}2<x_n<p_n-\prod_{k=1}^{n-1}p_k\bigg\}.
$$
Let $l_n$ be a positive integer such that
 $l_np_1\cdots p_{n-1}\le p_n/2<(l_n+1)p_1\cdots p_{n-1}$.
 We now set 
 $$
 Z_n:=\{x_n\in X_n \mid p_1\cdots p_{n-1}<x_n<p_n/2\}.
 $$
For each $n>0$, we define a probability measure $\mu_n$ on $X_n$ by the following conditions:
\roster
\item"---"
$\mu_n(j)=\mu_n(0)$ for each $j\le p_n/2$,
\item"---"
$\mu_n(j)=\mu_n(p_n-1)$ for each $j> p_n/2$,
\item"---"
$\frac{\mu_n(p_n-1)}{\mu_n(0)}=\lambda$.
\endroster
Of course, such a measure is defined uniquely. 
It is straightforward to verify that $H^2(\mu_n,\mu_n\circ T_n^{-1})\to 0$ as $n\to\infty$.
Passing to a countable subset in $P$, if necessary, we may assume  that the following three conditions are satisfied:
\roster
\item"$(\circ)$"
$\sum_{n=1}^\infty H^2(\mu_n,\mu_n\circ T_n^{-1})<\infty$,
\item"$(\bullet)$" $\mu_n(Y_n)>1-2^{-n-1}$ for each $n>0$,
\item"$(\star)$" $\mu_n(Z_n)>\frac 1{2(\lambda+1)}$  for each $n>0$.
\endroster
By Theorem~A, in view of $(\circ)$, $T$ is $\mu$-nonsingular\footnote{Though the topological system $(X,T)$ is a topological odometer, i.e. a minimal rotation on a monothetic compact  totally disconnected Abelian group, the nonsingular system $(X,\mu,T)$ should not be confused with the nonsingular product odometers which are well studied in the literature (see, e.g., \cite{HaOs}, \cite{DaSi}, \cite{Sc}) because $\mu$ does not split into infinite product when $X$ is written in the product form suitable for the odometer ``addition with carry".}, where $\mu=\bigotimes_{n=1}^\infty\mu_n$.

\proclaim{Proposition 2.12} $(X,\mu, T)$ is of Krieger's  type $III_\lambda$.
\endproclaim
\demo{Proof}
It follows from the Kolmogorov 0-1 law that $T$ is $\mu$-ergodic.
We claim that $T$ is of type $III_\lambda$.
Since $\log\frac{d\mu\circ T^{-1}}{d\mu}(x)\in\{n\log \lambda\mid n\in\Bbb Z\}$ at a.e. $x\in X$, it suffices to show that $\log\lambda$ is an essential value of the Radon-Nikodym cocycle $\rho_\mu$ of $T$.
Fix $n>0$ and take a Borel subset $B\subset X_1\times\cdots\times X_n$.
We now set $A:=B\times Z_{n+1}\times Y_{n+2}\times Y_{n+3}\times\cdots\subset X$.
Then $A$ is a Borel subset of the cylinder $[B]_{1}^n$.
Of course, $T^{p_1\cdots p_n}[B]_{1}^n=[B]_{1}^n$ and hence $T^{p_1\cdots p_nl_{n+1}}A\subset[B]_{1}^n$.
Since 
\roster
\item"---" $T_m^{p_1\cdots p_n}=I$ for each $m=1,\dots,n$,
\item"---"
$\frac{d\mu_{n+1}\circ T^{p_1\cdots p_nl_{n+1}}}{d\mu_{n+1}}(x_{n+1})=\lambda$ if $x_{n+1}\in Z_{n+1}$ and
\item"---"
$\frac{d\mu_m\circ T^k}{d\mu_m}(x_m)=1$ if $x_m\in Y_m$ and $0\le k\le p_1\cdots p_{m-1}$ and every $m>n+1$, 
\endroster
it follows that for each $x=(x_m)_{m=1}^\infty\in A$,
$$
\frac{d\mu\circ T^{p_1\cdots p_nl_{n+1}}}{d\mu}(x)=\prod_{m=1}^\infty\frac{d\mu_m\circ T_m^{p_1\cdots p_nl_{n+1}}}{d\mu_m}(x_m)
=\lambda.
$$
We also note that $\mu(A)>\frac{\mu([B]_1^n)}{2(\lambda+1)}\prod_{m=1}^\infty(1-2^{-m-1})$ in view of $(\bullet)$ and $(\star)$.
Hence, $\log\lambda$ is an essential value of $\rho_\mu$ by Lemma~D. 
\qed
\enddemo

\comment

\subhead 2.6. Calibration of the quasi-invariant measure for IDPFT systems
\endsubhead
Let $(X,\goth B,\mu,T)=\bigotimes_{n=1}^\infty(X_n,\goth B_n,\mu_n, T_n)$, $\mu_n\sim\nu_n$ and $\nu_n\circ T_n^{-1}=\nu_n$ for each $n>0$
and $\mu\circ T^{-1}\sim\mu$.
We
let
 $\phi_n:=\frac{d\mu_n}{d\nu_n}$ and
 $$
 t_{\sup}:=\sup\{t>0\mid \phi_n^t\in L^1(X_n,\nu_n)\text{ for all $n>0$}\}.
 $$
 It follows from the H{\"o}lder inequality that if $0<t_1<t_2$ and $\phi_n^{t_2}\in L^1(X_n,\nu_n)$
then $\phi_n^{t_1}\in L^1(X_n,\nu_n)$ and $\|\phi_n^{t_1}\|_1\le\|\phi_n^{t_2}\|_1^{\frac {t_1}{t_2}}$,
 where the norm is considered with respect to $\nu$.
 For $t\in(0, t_{\sup})$, define
  a probability measure $\mu_n^{(t)}$ on $X_n$
  by setting
  $$
\mu_n^{(t)}\sim\nu_n\qquad\text{and}\qquad\frac{d\mu_n^{(t)}}{d\nu_n}=\frac{\phi_n^t}{\|\phi_n^t\|_1}=:\phi_n^{(t)}.
$$
 We now define a probability measure $\mu^{(t)}$ on $X$ by setting
  $\mu^{(t)}:=\bigotimes_{n=1}^\infty\mu_n^{(t)}$.
  Let
$$
\text{Cal}(\mu):=\big\{t\in (0,t_{\sup})\mid\mu^{(t)}\circ T^{-1}\sim\mu^{(t)}\big\}.
$$
\definition{Definition 2.13}
If $t\in \text{Cal}(\mu)$ then we call $\mu^{(t)}$ the $t$-{\it calibration} of $\mu$.
\enddefinition

\proclaim{Proposition 2.14} 
\roster
\item"\rom{(i)}"  $\text{\rom{Cal}}(\mu)\supset\{t\in[1,t_{\sup})\mid \prod_{n=1}^\infty\|\phi_n^t\|_1<\infty\}$.
\item"\rom{(ii)}" For $t\in\text{\rom{Cal}}(\mu)$, we have that
$\mu^{(t)}\perp\mu$ if and only if $\prod_{n=1}^\infty\frac{\|\phi_n^{t/2}\|_1}{\|\phi_n^{t/2}\|_2}=0$.
\endroster
\endproclaim
\demo{Proof}
(i) Take $t\in[1,t_{\sup})$.
Applying  Jensen's inequality, we obtain that
$$
\align
\int_{X_n}\sqrt{\phi_n^{(t)}\cdot \phi_n^{(t)}\circ T_n^{-1}}\,d\nu_n&=
\frac 1{\|\phi_n^t\|_1}
\int_{X_n}\Big(\sqrt{\phi_n\cdot \phi_n\circ T_n^{-1}}\,\Big)^t\,d\nu_n\\
&\ge
 \frac 1{\|\phi_n^t\|_1}\bigg(\int_{X_n}\sqrt{\phi_n\cdot \phi_n\circ T_n^{-1}}\,d\nu_n\bigg)^t.
 \endalign
$$
Therefore,
$$
\prod_{n=1}^\infty\int_{X_n}\sqrt{\phi_n^{(t)}\cdot \phi_n^{(t)}\circ T_n^{-1}}\,d\nu_n\ge
\frac{\Big(\prod_{n=1}^\infty\int_{X_n}\sqrt{\phi_n\cdot \phi_n\circ T_n^{-1}}\,d\nu_n\Big)^t}{\prod_{n=1}^\infty\|\phi_n^t\|_1}>0
$$
in view of \thetag{2-2} and the fact that $\prod_{n=1}^\infty\|\phi_n^t\|_1<\infty$.
It remains to apply \thetag{2-2} again to conclude that $\mu^{(t)}\circ T\sim \mu^{(t)}$.

(ii) follows  directly from \thetag{2-3}.
\qed
\enddemo

\endcomment

\head 3. Gaussian dynamical systems\endhead

\subhead 3.1. Integration in Hilbert spaces
\endsubhead
Let $\Cal H$ denote a separable infinite dimensional real Hilbert space. Given a Borel probability measure $\mu$ on $\Cal H$, we denote by $\widehat\mu$ the {\it characteristic functional} of $\mu$, i.e.
$$
\widehat\mu(y):=\int_\Cal He^{i\langle x,y\rangle}d\mu(x), \quad y\in\Cal H.
$$
We note that each Borel probability measure on $\Cal H$ is defined completely by its characteristic functional.
If there is a vector $h\in\Cal H$ and a bounded linear operator $B>0$ in $\Cal H$ such that $\widehat\mu(y)=e^{i\langle h,y\rangle-\frac12\langle By,y\rangle}$ for all $y\in\Cal H$ then $\mu$ is called  the {\it (non-degenerated) Gaussian measure with covariance operator $B$ and mean $h$}.
Then for each $t\in\Bbb R$ and $y\in \Cal H$,
$$
\int_{\Bbb R}e^{its}\,d(\mu\circ\langle\cdot,y\rangle^{-1})(s)=
\int_{\Cal H}e^{i\langle x,ty\rangle}\,d\mu(x)=e^{it\langle h,y\rangle-\frac12t^2\langle By,y\rangle}.
$$
Therefore the continuous linear functional $\Cal H\ni x\mapsto\langle x,y\rangle$ has normal distribution $\Cal N_{\langle h, y\rangle,\langle By,y\rangle}$.
In particular, each continuous linear functional belong to $L^2(\Cal H,\mu_B)$.
By the Minlos-Sazonov theorem, $B$ is a nuclear operator, i.e. tr$(B)<\infty$ \cite{Sk}.
Conversely, each strictly positive nuclear operator $B$ in $\Cal H$ determines a  unique Gaussian measure on $\Cal H$ with zero mean and covariance operator $B$.
We denote this measure by $\mu_B$.
Thus $\widehat{\mu_B}(y)= e^{-\frac12\langle By,y\rangle}$ for all $y\in\Cal H$.
We note that  $\int_{\Cal{H}}\langle x,y\rangle d\mu_B(y)=0$ for each  $h\in\Cal{H}$.
It is well known that
$$
\int_\Cal H\langle x,y\rangle\langle z,y\rangle\, d\mu_B(y)=\langle Bx,z\rangle
$$
and  hence tr$(B)=\int_\Cal H\|y\|^2\,d\mu_B(y)$.
We now let $\Cal H_0:=B^{\frac12}\Cal H\subset \Cal H$ and define an inner product and the corresponding  norm on 
$\Cal H_0$ by setting
$$
\langle x,y\rangle_0:=\langle B^{-\frac12}x,B^{-\frac12}y\rangle\text{ and $\|x\|_0^2:=\langle x,x\rangle_0$ \  for } x,y\in\Cal H_0.
$$
Then $(\Cal H_0,\langle .,.\rangle_0)$ is a Hilbert space.
We now show that there is  a canonical isometric embedding of $\Cal H_0$ into $L^2(\Cal H,\mu_B)$.
For that, we first take  $\theta\in B\Cal H\subset\Cal H_0$.
Then  the mapping
$$
l_\theta:\Cal H\ni y\mapsto\langle B^{-1}\theta,y\rangle
$$
is a continuous linear functional on $\Cal H$.
Moreover, for all $\theta,\eta\in B\Cal H$,
$$
\langle l_\theta,l_\eta\rangle_{L^2(\Cal H,\mu_B)}=\int_\Cal H
\langle B^{-1}\theta,y\rangle\langle B^{-1}\eta,y\rangle\,d\mu_B(y)=\langle\theta, B^{-1}\eta
\rangle=
\langle\theta, \eta
\rangle_0.
$$
In particular,  the linear mapping
$$
l:B\Cal H\ni\theta\mapsto l_\theta\in L^2(\Cal H,\mu_B)
$$
is isometric\footnote{If $B\Cal H$ is furnished with $\|.\|_0$.}.
We note that $B\Cal H$ is dense in $\Cal H_0$.
Indeed,  since the linear span $\Cal L$ of the orthonormal basis in $\Cal H$ consisting
of eigenvectors for $B$ is dense in $\Cal H$, it follows that $B^{\frac12}\Cal L$ is dense
in $\Cal H_0$ because  $B^{\frac12}$ is an isometric isomorphism of $(\Cal H,\langle,.\rangle)$
onto $(\Cal H_0,\langle,.\rangle_0)$.
It remains to observe that $B^{\frac12}\Cal L=B\Cal L=\Cal L$.
Since $B\Cal H$ is dense in $\Cal H_0$,  the isometry $l$ extends by continuity to an isometry from $\Cal H_0$ to $L^2(\Cal H,\mu_B)$.
Thus, for each $y\in\Cal H_0$, there is a sequence $(\theta_n)_{n=1}^\infty$ of elements from $B\Cal H$ such that $\|y-\theta_n\|_0\to 0$ and the sequence $(l_{\theta_n})_{n=1}^\infty$ converges to some element $l_y\in L^2(\Cal H,\mu_B)$.
Hence a subsequence $(l_{\theta_{n_k}})_{k=1}^\infty$ converges  to $l_y$ almost everywhere.
Let $D_y$ denote the set of all $x\in\Cal H$ such that the sequence $(l_{\theta_{n_k}}(x))_{k=1}^\infty$ converges. 
It is easy to verify that $D_y$ is a (Borel) linear subspace of $\Cal H$,  $\mu_B(D_y)=1$ and $l_y$ is linear on $D_y$.
That is why $l_y$  is often called a {\it  measurable linear functional on $\Cal H$}.
Moreover, $D_y\supset\Cal H_0$ and $l_y$ is defined uniquely by the restriction to $\Cal H_0$ though $\mu_B(\Cal H_0)=0$.
It is often convenient to write $\langle B^{-1}y,x\rangle$ instead of  $l_y(x)$ for $\mu_B$-a.a.  $x\in \Cal H$.
We note that the distribution of $l_\theta$ is $\Cal N_{0,\langle\theta,B^{-1}\theta\rangle}=\Cal N_{0,\|\theta\|_0^2}$ for each $\theta\in B\Cal H$.
Passing to a limit we obtain that the distribution of $l_y$ is $\Cal N_{0,\|y\|_0^2}$ for each $y\in\Cal H_0$.

For each $y\in \Cal H$, we denote by $L_y$ the rotation by $y$, i.e. $L_yx=x+y$ for all $x\in\Cal H$. 
By the Cameron-Martin theorem (see \cite{Gu, Corollary 7.4}, \cite{Sk}),  
$$
\gather
\Cal H_0=\{y\in\Cal H\mid \mu_B\sim\mu_B\circ L_y^{-1}\}\quad\text{and for each $y\in \Cal H_0$,}\\
\frac{d\mu_B\circ L_y^{-1}}{d\mu_B}(x)=
e^{\langle B^{-1}y,x\rangle-\frac12\|y\|_0^2}
\quad\text{at a.e. $x\in\Cal H$}. \tag3-1
\endgather
$$

\subhead 3.2. Fock space and exponential map\endsubhead
Given a separable Hilbert space $\Cal K$,
{\it the (bosonic) Fock space $\Cal F(\Cal K)$ built over $\Cal K$} is the Hilbert space
$\bigoplus_{n=0}^\infty\Cal K^{\odot n}$.
The subspace $\Cal K^{\odot n}$ of $\Cal F(\Cal K)$ is called {\it the $n$-chaos} in $\Cal F(\Cal K)$, $n\in\Bbb Z_+$.
Given $h\in\Cal K$, we let $\exp_h:=\bigoplus_{n=0}^\infty\frac{h^{\otimes n}}{\sqrt{n!}}\in\Cal F(\Cal K)$.
In particular, $\exp_0=(1,0,0,\dots,)$ is called {\it the vacuum vector} in $\Cal F(\Cal H)$.
The map $\exp:\Cal K\ni h\mapsto \exp_h\in\Cal F(\Cal K)$ is called {\it the exponential map}.
It satisfies the following \cite{Gu}:
\roster
\item"(i)"
$\exp$ is continuous,
\item"(ii)" $\langle \text{exp}_h,\text{exp}_k\rangle_{\Cal F(\Cal K)}=e^{\langle h,k\rangle_{\Cal K}}$ for all $h,k\in\Cal K$,
\item"(iii)"  the set $\{\text{exp}_h\mid h\in\Cal K\}$ is  linearly independent and total
in $\Cal F(\Cal K)$.
\endroster
Given an orthogonal operator $V$ in $\Cal K$, we can define a linear operator $\exp V$ of $\Cal F(\Cal K)$, called {\it the second quantization of $V$}, by setting
$$
(\exp V)  h^{\otimes n}:= (Vh)^{\otimes n} \quad \text{for all $n\ge 0$ and $h\in\Cal K$.}
$$
Then $\exp V$ preserves each chaos in $\Cal F(\Cal K)$ and the restriction of $\exp V$ to the first chaos is $V$.
Of course, $(\exp V)\exp_h=\exp_{Vh}$ for each $h\in\Cal K$.
The most important property of the Fock spaces is the following one: given a decomposition $\Cal K=\bigoplus_{j=1}^\infty\Cal K_j$ of $\Cal K$ into an orthogonal sum of subspaces $\Cal K_j$, there is a unique unitary isomorphism $\Phi$ of $(\Cal F(\Cal K),\exp_0)$ onto
$\bigotimes_{j=1}^\infty(\Cal F(K_j),\exp_0)$ such that
$$
\Phi(\exp_{\bigoplus_{j=1}^\infty h_j})=\bigotimes_{j=1}^\infty\Phi(\exp_{h_j})
$$
for each vector $\bigoplus_{j=1}^\infty h_j\in\Cal K$ such that $h_j=0$ for all but finitely many $j$ \cite{Gu, Proposition~2.3}.

Denote  the orthogonal group  of $\Cal K$ by $\Cal O(\Cal K)$.
Let $ \text{Aff}(\Cal K):=\Cal K\rtimes \Cal O(\Cal K)$ stand for   the group of  affine operators in $\Cal K$.
We recall that 
an operator
 $A= (f,V)\in \text{Aff}(\Cal K)$ acts on $\Cal K$
  by the formula
  $Ah:=f+V h$.
   One can verify that the multiplication law in Aff$(\Cal K)$ is given by:
   $$
   (f,V)(f',V') := (f+V f',V V').
   $$
   We note that Aff$(\Cal K)$ is a Polish group if endowed with the product of the norm topology on $\Cal K$  and the weak operator topology on 
  $ \Cal O(\Cal K)$.
We  recall  the well-known {\it Weyl unitary representation} $W=(W_{(f,V)})_{(f,V)\in\text{Aff}(\Cal K)}$ of Aff$(\Cal K)$
in  $\Cal F(\Cal K)$ \cite{Gu, \S2.2}:
$$
W_{(f,V)}\exp_h:=e^{-\langle f,Vh\rangle_\Cal K-\frac12\|f\|^2_\Cal K}\exp_{f+Vh},\qquad  h\in\Cal K.\tag3-2
$$
It is well defined due to (ii) and (iii).
Of course, $W_{(0,V)}=\exp V$ for each $V\in\Cal O(\Cal K)$.

By \cite{Gu, Theorem~7.1}, there is a unique (canonical) unitary isomorphism of $L^2(\Cal H,\mu_B)$ with $\Cal F(\Cal H_0)$ such that\footnote{For simplicity sake, we will write $L^2(\Cal H,\mu_B)=\Cal F(\Cal H_0)$ and hence identify $\exp_h$ with an
 $L^2$-function on $(\Cal H,\mu_B)$, $h\in\Cal H_0$.}
 $$
 \text{exp}_h(x):=e^{\langle B^{-1}h,x\rangle-\frac12\|h\|_0^2},\quad \text{for a.e. $x\in\Cal H$.}\tag3-3
 $$
Moreover, the map $\Cal H_0\ni h\mapsto l_h\in L^2(\Cal H,\mu_B)$ identifies (isometrically)
$\Cal H_0$ with the first chaos in $L^2(\Cal H,\mu_B)$.
It follows from \thetag{3-1}  and \thetag{3-3} that
 $$ \text{exp}_h=\frac{d\mu_B\circ L_h^{-1}}{d\mu_B}\quad\text{ for each $h\in\Cal H_0$.}
 \tag3-4
 $$
It is straightforward to verify that the following additional properties for $\exp$ hold:
\roster
\item"(iv)" exp$_h>0$ for each $h\in \Cal H_0$,
\item"(v)"  exp$_h\in\bigcap_{p=1}^\infty L^p(\Cal H,\mu_B)$
because  the map $\Cal H\ni x\mapsto\langle B^{-1}h,x\rangle -\frac12\|h\|^2_0$ has normal distribution
$\Cal N_{-\frac12\|h\|_0^2,\|h\|_0^2}$ and \thetag{3-3} holds,
 \item"(vi)"$\|\text{exp}_h\|_1=1$ for each $h\in \Cal H_0$,
\item"(vii)"
the cone $\{\sum_{k=1}^na_k\text{exp}_{h_k}\mid a_1,\dots,a_n>0, h_1,\dots,h_n\in\Cal H_0, n\in\Bbb N\}$ is dense in the cone $L^2_+(\Cal H,\mu_B)$ of non-negative functions from $L^2(\Cal H,\mu_B)$,
\item"(viii)"
$\text{exp}_h\cdot \text{exp}_k=e^{\langle h,k\rangle_0}\text{exp}_{h+k}$
for all $h,k\in\Cal H_0$ and hence
\item"{(ix)}"
$\sqrt{\text{exp}_h}=e^{-\frac 18\|h\|_0^2}\text{exp}_{h/2}
$
for each $h\in\Cal H_0$,
\item"{(x)}" $\exp_h\circ L_f^{-1}=e^{-\langle B^{-1}h,f\rangle}\exp_h=e^{-\langle h ,f\rangle_0}\exp_h$ for all $h,f\in\Cal H_0$.
\endroster

\remark{Remark \rom{3.1}}\roster\item"(i)" We recall that $\Cal H_0$ is determined by the pair $(\Cal H,B)$ (see \S3.1).
Conversely,   if $\Cal H_0$ is given beforehand  as an abstract  Hilbert space, then it 
determines uniquely the probability space $(\Cal H,\mu_B)$  for some pair $(\Cal H, B)$ such that $\Cal H_0=B^{\frac12} \Cal H$.
Indeed, if there is another Hilbert space $\Cal K$ and a non-degenerated nuclear operator $C>0$ on $\Cal K$ such that the  space $\Cal K_0:=C^{\frac 12}\Cal K$ furnished with the corresponding Hilbert norm is unitarily  isomorphic to $\Cal H_0$ via some unitary isomorphism $\Psi$ then according to 
\cite{Gu, Theorem~7.1} and \thetag{3-3}, there is a unique unitary isomorphism  $\Phi$ of  $L^2(\Cal H,\mu_B)$ with $L^2(\Cal K,\mu_C)$ which  maps  $\exp_h$ onto  $\exp_{\Psi^{-1} h}$ for each $h\in\Cal H_0$.
Hence in view of (vii), $\Phi$ maps  $L^2_+(\Cal K,\mu_C)$ onto $L^2_+(\Cal K,\mu_C)$.
Moreover, $\Phi 1=1$.
Therefore $\Phi$ is {\it spacial}, i.e. there is a measure preserving isomorphism $\theta:(\Cal H,\mu_B)\to(\Cal K,\mu_C)$ such that $\Phi h=h\circ\theta^{-1}$ for each $h\in\Cal H$.
\item"(ii)"
Another useful observation is that  given a Hilbert space $\Cal K_0$, there is another Hilbert space $\Cal K\supset \Cal K_0$ and a  nuclear operator $C$ of $\Cal K$ such that  $C^{\frac12}$ is a unitary isomorphism of $\Cal K$ onto $\Cal K_0$.
\endroster
\endremark

\remark{Remark \rom{3.2}} Given  a decomposition $\Cal H_0=\bigoplus_{j=1}^\infty\Cal H_{0,j}$ of $\Cal H_0$ into an orthogonal sum of subspaces $\Cal H_{0,j}$, consider
the corresponding decomposition $\Cal H=\bigoplus_{j=1}^\infty\Cal H_{j}$ of $\Cal H$ into an orthogonal sum of subspaces $\Cal H_j:=B^{-\frac 12}\Cal H_{j,0}$, $j\in\Bbb N$.
Let $P_j:\Cal H\to\Cal H_j$ denote the orthogonal projection of $\Cal H$ onto $\Cal H_j$ and let $B_j:=P_jBP_j^*$.
Then $B_j:H_j\to H_j$ is a nuclear operator and $B_j^{\frac12}\Cal H_j=\Cal H_{0,j}$ for each $j\in\Bbb N$.
Moreover,  $(\Cal H,\mu_B)$ splits into the  direct product $(\Cal H,\mu_B)=\bigotimes_{j=1}^\infty(\Cal H_j,\mu_{B_j})$ of Gaussian probability spaces $(\Cal H_j,\mu_{B_j})$ in such a way that
$\{\exp_h\mid h\in\Cal H_{0,j}\}$ is total in $L^2(H_j,\mu_{B_j})$ and $\mu_{B_j}=\mu_B\circ P_j^*$
for each $j\in\Bbb N$.
\endremark

\subhead 3.3. Nonsingular Gaussian action of \rom{Aff}$\,\Cal H_0$
\endsubhead
Let $(Y,\goth C,\nu)$ be a standard nonatomic probability space.
Denote by $\Cal U(L^2(Y,\nu))$ the group of unitary operators in $L^2(Y,\nu)$ and 
by $\Cal U_\Bbb R(L^2(Y,\nu))$ the subgroup of unitaries  that preserve  the subspace $L^2_{\Bbb R}(Y,\nu)$ of real valued functions in $L^2(Y,\nu)$.
Let
$$
U:\text{Aut}(Y,\nu)\ni T\mapsto U_T\in\Cal U_\Bbb R(L^2(Y,\nu))
$$
stand for the unitary Koopman representation of $\text{Aut}(Y,\nu)$ in $L^2(Y,\nu)$.
We recall that $U_Tf:=f\circ T^{-1}\sqrt{\frac{d\mu\circ T^{-1}}{d\mu}}$ for all $f\in L^2(Y,\nu)$.
 The following results are well known:
 \roster
 \item"($\bullet)$"  $\{U_T\mid T\in \text{Aut}(Y,\nu)\}=\{V\in \Cal U_\Bbb R(L^2(Y,\nu))\mid VL^2_+(Y,\nu)=L^2_+(Y,\nu)\}$.
  \item"($\circ)$"  $\{U_T\mid T\in \text{Aut}_0(Y,\nu)\}=\{V\in \Cal U_\Bbb R(L^2(Y,\nu))\mid VL^2_+(Y,\nu)=L^2_+(Y,\nu), V1=1\}$.
 \endroster
 We also note that $U$ is one-to-one and the image of $U$ is closed in $\Cal U_\Bbb R(L^2(Y,\nu))$ in the weak (and the strong) operator topology.

  Let $\Bbb R^*$ denote the multiplicative group of reals.
  It is straightforward to verify that for each $t\in\Bbb R^*$, the map $\alpha_t:\text{Aff}(\Cal H_0)\to\text{Aff}(\Cal H_0)$ given by
  $$
(f,V)\mapsto  \alpha_t(f,V):=(tf,V)\tag3-5
  $$
 is a continuous automorphism of $\text{Aff}(\Cal H_0)$.
 Moreover, $\alpha_{t_1}\alpha_{t_2}=\alpha_{t_1t_2}$ for all $t_1,t_2\in\Bbb R^*$.

It is straightforward to verify that for each $A\in\text{Aff}(\Cal H_0)$, the corresponding Weyl unitary operator $W_{A}$ (see \thetag{3-2}) preserves the cone 
$$
\bigg\{\sum_{k=1}^na_k\,\text{exp}_{h_k}\,\Big|\, a_k>0, h_k\in\Cal H_0,\text{for  each $k=1,\dots,n$ and } n\in\Bbb N\bigg\}.
$$
Hence it preserves $L^2_+(\Cal H,\mu_B)$ in view of (vii) from \S3.2.
Therefore by $(\bullet)$, there is a (unique) transformation $T_{A}\in\text{Aut}(\Cal H,\mu_B)$ such that
$U_{T_A}=W_{\alpha_{1/2}(A)}$.

\definition{Definition 3.3} $T_{A}$ is called  the {\it nonsingular Gaussian transformation generated by} $A\in \text{Aff}(\Cal H_0)$. 
\enddefinition

Since the image of $\text{Aff}(\Cal H_0)$ under the unitary Weyl representation is closed in the unitary group of the space $L^2(\Cal H,\mu_B)$ \cite{Gu, Theorem~2.1},
it follows  that  the group $\{T_A\mid A\in \text{Aff}(\Cal H_0)\}$ of nonsingular Gaussian transformations is closed in Aut$(\Cal H,\mu_B)$.

\proclaim{Proposition 3.4}
\roster
\item "\rom{(i)}"If $V\in \Cal O(\Cal H_0)$ then  $T_{(0,V)}$ is the usual (classic) measure preserving Gaussian transformation generated by the orthogonal operator $V$, i.e.
$U_{T_{(0,V)}}=\exp V$ (see \cite{LePaTh, Lemma~2}).
\item"\rom{(ii)}"
 If $f\in\Cal H_0$ then $T_{(f,I)}=L_f$.
\endroster
\endproclaim

\demo{Proof}
(i) We note that 
$$
U_{T_{(0,V)}}\text{\rom{exp}}_h=W_{(0,V)}\text{\rom{exp}}_h=\text{exp}_{Vh}=(\exp V)\exp_h.
$$
Hence $U_{T_{(0,V)}}=\exp V$.
%If $x\in \Cal H_0$ then
%$$
%\text{exp}_{Vh}(x)=e^{\langle Vh,x\rangle_0-\frac12\|Vh\|_0^2}=e^{\langle h,V^{-1}x\rangle_0-\frac12\|h\|_0^2}=\text{exp}_{h}(V^{-1}x).
%$$

(ii)
Using \thetag{3-4} and (viii)--(x) from \S3.2  we obtain that
$$
\align
U_{L_f}\text{exp}_h&=\sqrt{\frac{d\mu_B\circ L_f^{-1}}{d\mu_B}}\,\text{exp}_h\circ L_f^{-1}\\
&=
\sqrt{\text{exp}_{f}} \, e^{- \langle f,h\rangle_0}\exp_h\\
%e^{\langle B^{-1}h,x-f\rangle-\frac12\|h\|_0^2}\\
&=
e^{-\frac 18\|f\|_0^2}\text{exp}_{f/2}\, e^{- \langle h,f\rangle_0}\exp_h\\
&=e^{-\frac 18\|f\|_0^2-\langle f,h\rangle_0+\frac12\langle f,h\rangle_0}\text{exp}_{f/2+h}.
\endalign
$$
Hence $U_{L_f}\text{exp}_h=e^{-\frac 18\|f\|_0^2-\frac12\langle f,h\rangle_0}\text{exp}_{f/2+h}=W_{(f/2,I)}\text{exp}_h=U_{T_{(f,I)}}\text{exp}_h$.
It follows that $T_{(f,I)}=L_f$.
\qed
\enddemo

\proclaim{Corollary 3.5} Every nonsingular Gaussian transformation $T_{(f,V)}$ is the composition of the classic $\mu_B$-preserving Gaussian transformation  $T_{(0,V)}$ and
a $\mu_B$-nonsingular translation $L_f=T_{(f,I)}$ which is totally dissipative.\footnote{Let $\Cal K$ stand for the orthogonal complement in $\Cal H$ to the 1-dimensional subspace generated by $f$.  Then the set $\{sf+k\mid 0\le s<1, \, k\in\Cal K\}\subset\Cal H$ is a Borel fundamental domain for $L_f$.}
These two transformations  commute if and only if  $Vf=f$.
\endproclaim

\remark{Remark \rom{3.6}} Let $(X,\goth B,\mu)$ be a standard $\sigma$-finite nonatomic measure space.
Let a transformation $S\in\text{Aut}(X,\mu)$ be such that $\sqrt{\frac{d\mu\circ S^{-1}}{d\mu}}-1\in L^2(X,\mu)$.
Then a nonsingular Poisson suspension $S_*$ of $S$ is well defined on a standard probability space $(X^*,\goth B^*,\mu^*)$ \cite{DaKoRo1}.
Let $A:=(U_S, \sqrt{\frac{d\mu\circ S^{-1}}{d\mu}}-1))\in\text{Aff}(L^2(X,\mu))$.
It was shown in  \cite{DaKoRo1} that $U_{S_*}$ is unitarily equivalent to $W_A$.
It follows that  each nonsingular Poisson transformation is unitarily equivalent to a nonsingular Gaussian transformation: $S_*$  is unitarily equivalent  to  $T_{\alpha_{2}(A)}$ (see \thetag{3-5}).
We do not know if the converse is true even in the classic (finite measure preserving) case.

\endremark

It is well known that the transformation group $\{T_{(f,I)}\mid f\in\Cal H_0\}\subset\text{Aut}(\Cal H,\mu_B)$ is ergodic (see \cite{Gu}, \cite{Sk}).
However  Krieger's type of it has not been determined so far.
We will show that it is  type $III_1$, i.e. a  dense countable subgroup of it is of type $III_1$ (hence every dense countable subgroup is of type $III_1$).

\proclaim{Theorem 3.7} $\{T_{(f,I)}\mid f\in\Cal H_0\}$  is of type $III_1$.
\endproclaim
\demo{Proof}
Let $\{e_n\mid n\in\Bbb N\}$ be an orthonormal basis of $\Cal H$ consisting of the eigenvectors of $B$. 
Then $Be_n=\lambda_ne_n$, $\lambda_n>0$ for each $n\in\Bbb N$ and $\sum_{n=1}^\infty\lambda_n<\infty$.
Denote by $\Gamma$ the group generated by translations $L_{\sqrt{\lambda_k}e_k}$ for all  $k\in\Bbb N$.
Then $\Gamma$ is an  ergodic countable Abelian subgroup of $\text{Aut}(\Cal H,\mu_B)$.
We will show that $\Gamma$ is of type $III_1$.

Denote by $\goth B_n$ the smallest Borel $\sigma$-algebra on $\Cal H$ such that the map $\Cal H\ni x\mapsto\langle x, e_k\rangle\in \Bbb R$ is $\goth B_n$-measurable for each $k=1,\dots,n$.
Then $\goth B_1\subset\goth B_2\subset\cdots$ and the union $\bigcup_{n>0}\goth B_n$ is dense in $\goth B$.
We deduce from \thetag{3-3} and \thetag{3-4} that  for each $n>0$,
$$
\log\frac{d\mu_B\circ L_{\sqrt{\lambda_{n+1}}e_{n+1}}}{d\mu_B}(x)=\frac{\langle x, e_{n+1}\rangle}{\sqrt{\lambda_{n+1}}}-\frac12.
$$
Take $a\in\Bbb R$ and $\epsilon>0$.
We now let
$$
D_n:=\bigg\{x\in\Cal H\,\bigg|\,  a+\frac12 -\epsilon <\frac{\langle x, e_{n+1}\rangle}{\sqrt{\lambda_{n+1}}}  <a+\frac12+\epsilon \bigg\}.
$$
Since $e_{n+1}\perp e_k$ for each $k=1,\dots,n$ and  the random variable $\Cal H\ni x\mapsto \langle x,e_k\rangle\in\Bbb R$ is Gaussian for all $k=1,\dots,n+1$ (and the joint distribution is distributions are also Gaussian), it follows that
 $D_n$ is independent of $\goth B_n$.
Moreover, the measure
$$
\align
\mu_B(D_n)&=\frac1{\sqrt{2\pi\lambda_{n+1}}}\int_{(a+\frac12-\epsilon,a+\frac12+\epsilon)\cdot\sqrt{\lambda_{n+1}}}e^{-\frac{t^2}{2\lambda_{n+1}}}dt\\
&=
\frac1{\sqrt{2\pi}}\int_{(a+\frac12-\epsilon,a+\frac12+\epsilon)}e^{-\frac{t^2}2}dt
\endalign
$$
of $D_n$ does not depend on $n$.\footnote{We use here the fact that the random variable $\langle\cdot,e_{n+1}\rangle$ has normal distribution $\Cal N_{0, \lambda_{n+1}}$.}
We denote it by $\delta>0$.
Then for each subset $A\in\goth B_n$, we have that $L_{e_{n+1}}A=A$ and hence
\roster
\item
"---" $(A\cap D_n)\cup L_{e_{n+1}}(A\cap D_n)\subset A$, 
\item"---" $\mu_B(A\cap D_n)=\mu_B(A)\mu_B(D_n)=\delta\mu_B(A)$
and 
\item"---"
$\log\frac{d\mu_B\circ L_{e_{n+1}}}{d\mu_B}(x)=a\pm\epsilon$ for each $x\in A\cap D_n$.
\endroster
It follows from Lemma D (see \S2.2) that $a$ is an essential value of the logarithm of the Radon-Nikodym cocycle of $\Gamma$.
Since $a$ is an arbitrary element of $\Bbb R$, the Radon-Nikodym cocycle is ergodic, i.e.
$\Gamma$ is of type $III_1$.
\qed
\enddemo

\subhead 3.4. When nonsingular Gaussian systems are of type $II_1$
\endsubhead
We recall a standard definition.

\definition{Definition 3.8}
Given $V\in \Cal O(\Cal H_0)$, we say that a vector $f\in\Cal H_0$ is a {\it $V$-coboundary} if there is $a\in\Cal H_0$ such that $f=a-Va$.
\enddefinition
In this subsection we prove the following
statement (cf. \cite{DaKoRo1, Proposition~6.4} and \cite{ArIsMa}).

\proclaim{Theorem 3.9} Let $(f,V)\in\text{\rom{Aff}}(\Cal H_0)$.
For $n\in\Bbb Z$, we define $f^{(n)}\in\Cal H_0$ by setting $(f,V)^n=(f^{(n)}, V^n)$.
The following are equivalent:
\roster
\item"\rom(i)" $T_{(f,V)}$ admits an equivalent invariant probability measure.
\item"\rom(ii)" $f$ is a $V$-coboundary.
\item"\rom(iii)" The affine operator $(f,V)$ has a fixed point.
\item"\rom(iv)" The sequence $(f^{(n)})_{n\in\Bbb Z}$ is bounded in $\Cal H_0$.
\endroster
\endproclaim

\demo{Proof}  (ii)$\iff$(iv) is classic, see  \cite{BeKaVa, Proposition 2.2.9}, for a proof.

(ii)$\iff$(iii) is obvious because the equality $(f,V)a=a$ for some $a\in\Cal H_0$ means
$f+Va=a$, i.e. $f$ is a $V$-coboundary.

(ii)$\Longrightarrow$(i) 
In view of Proposition~3.2,
$$
\frac{d\mu_B\circ T_{(f,V)}^{-1}}{d\mu_B}=\frac{d(\mu_B\circ T_{(0,V)}^{-1})\circ T_{(f,I)}^{-1}}{d\mu_B}
=\frac{d\mu_B\circ T_{(f,I)}^{-1}}{d\mu_B}=\frac{d\mu_B\circ L_f^{-1}}{d\mu_B}.
$$
Therefore, by \thetag{3-4}, we obtain that
$$
\frac{d\mu_B\circ T_{(f,V)}^{-1}}{d\mu_B}=\text{exp}_f.
$$
Let $f=a-Va$ for some $a\in\Cal H_0$.
We claim that 
$$
\text{exp}_f=\frac{\text{exp}_a}{\text{exp}_a\circ T_{(f,V)}^{-1}}.\tag3-6
$$
Indeed, applying  Proposition~3.4 and (viii) and (x) from \S3.2, we obtain that
$$
\align
\text{exp}_f\,\text{exp}_a\circ T_{(f,V)}^{-1}&=
\text{exp}_f\,\text{exp}_a\circ T_{(0,V)^{-1}}
\circ  L_f^{-1}\\
&=
\text{exp}_f\,((\exp V)\exp_a)\circ L_f^{-1}\\
&=
\exp_f\,\exp_{Va}\,e^{-\langle Va,f\rangle}\\
&=\text{exp}_{f+Va}\\
&=\text{exp}_{a}.
\endalign
$$
Since exp$_a\in L^1(\Cal H,\mu_B)$,  (i) follows from \thetag{3-6}.

(i)$\Longrightarrow$(iv)
We first note that for each $h\in\Cal H_0$,
$$
\|\sqrt{\text{exp}_h}\|_1
=e^{-\frac{\|h\|_0^2}{8}}.
$$
We now have
$$
\align
\langle (U_{T_{(f,V)}})^n1,1\rangle&=\Big\langle U_{T_{(f^{(n)},V^n)}}1,1\Big\rangle\\
&=\Bigg\langle\sqrt{\frac{d\mu_B\circ T_{(f^{(n)},V^n)}^{-1}}{d\mu_B}},1\Bigg\rangle\\
&=
\|\sqrt{\text{exp}_{f^{(n)}}}\|_1\\
&=e^{-\frac{\|f^{(n)}\|_0^2}{8}}.
\endalign
$$
Suppose that the sequence $(f^{(n)})_{n=1}^\infty$ is unbounded.
Then there is an increasing sequence $n_1<n_2<\cdots$ such that $\|f^{(n_k)}\|^2_0\to+\infty$ as $k\to\infty$.
Hence $\langle (U_{T_{(f,V)}})^{n_k}1,1\rangle\to 0$ as $k\to\infty$.
Since the operator $U_{T_{(f,V)}}$ is positive with respect to the cone $L^2_+(\Cal H,\mu_B)$, it
 follows that $U_{T_{(f,V)}}^{n_k}\to 0$ weakly as $k\to\infty$.
Since $T_{(f,V)}$ admits an equivalent invariant probability measure, $U_{T_{(f,V)}}$ is unitarily equivalent to the Koopman operator of a probability preserving transformation.
The latter does not have subsequences weakly converging to zero because $1$ is a fixed point of this operator. \qed
\enddemo

\remark{Remark 3.10} In fact, we  showed more: if  $f=a-Va$ and $\nu$ is a $\mu_B$-equivalent  $T_{(f,V)}$-invariant measure  then
$\frac{d\nu}{d\mu_B}=\exp_a$.
\endremark

\subhead 3.5. Gaussian transformations as IDPFT systems
\endsubhead
%Suppose that we have a decomposition $\Cal H_0=\bigoplus_{r=1}^\infty\Cal H_{0,r}$ of $\Cal H_0$ into a sum of mutually orthogonal closed subspaces $\Cal H_{0,r}$.
%One of the most important properties of the Fock spaces is that $\Cal F(\Cal H_0)$ ``splits'' canonically into the corresponding tensor product $\Cal F(\Cal H_0)=\bigotimes_{r=1}^\infty\Cal F(\Cal H_{0,r})$ in such  a way that $\Cal F(\Cal H_{0,r})$
%is generated by the total family $\{\text{exp}_h\mid h\in \Cal H_{0,r}\}$ for each $r\in\Bbb N$.
Suppose that we are given an affine operator $(f,V)\in\text{Aff}(\Cal H_0)$.
Suppose also that $V$ has no non-trivial invariant vectors.\footnote{Equivalently, the measure of maximal spectral type of $V$ has no atom at $1$.}
Using the spectral decomposition of $V$, as in the proof of Theorem~1.6, we can choose
 an orthogonal decomposition $\Cal H_0=\bigoplus_{r=1}^\infty\Cal H_{0,r}$ of $\Cal H_0$
 in such a way that  $V\Cal H_{0,r}=\Cal H_{0,r}$ and the orthogonal projection $f_r$ of $f$ onto $\Cal H_{0,r}$ is a $V$-coboundary for each $r\in\Bbb N$.
Let $V_r:=V\restriction\Cal H_r$.
Then $(f_r, V_r)\in\text{Aff}(\Cal H_{0,r})$ for each $r\in\Bbb N$
 and $(f,V)=\bigoplus_{r=1}^\infty (f_r,V_r)$.
Let $\Cal H_r$ and $\mu_r$ stand for the Hilbert space and a Gaussian measure on $\Cal H_r$ respectively such that $\Cal F(\Cal H_{0,r})$ is canonically isomorphic to
$L^2(\Cal H_r,\mu_r)$ (see Remark~3.1(ii)).
Then the standard probability space $(\Cal H,\mu_B)$ is isomorphic to the infinite product
$\bigotimes_{r=1}^\infty(\Cal H_r,\mu_r)$ according to Remark~3.2.
It follows that
$$
(\Cal H,\mu_B,T_{(f,V)})=\bigotimes_{r=1}^\infty(\Cal H_r,\mu_r, T_{(f_r,V_r)}).\tag3-7
$$
  Since $f_r$ is a $V_r$-coboundary, 
 there is $a_r\in\Cal H_{0,r}$ such that $f_r=a_r-V_ra_r$
 for each $r\in\Bbb N$.
By Theorem~3.9, the system $(\Cal H_r,\mu_r, T_{(f_r,V_r)})$ admits an equivalent invariant probability measure $\nu_r$.
Moreover, $\frac{d\mu_r}{d\nu_r}=\text{exp}_{-a_r}$ for each $r\in\Bbb N$ in view of Remark~3.10.
Thus, we have shown that each nonsingular Gaussian dynamical system $(\Cal H,\mu_B,T_{(f,V)})$ such that $V$ has no non-trivial invariant vectors
is  IDPFT (see \thetag{3-7}).
Therefore, Corollary~2.7 yields the following.

\proclaim{Corollary 3.11} If $V$ has no nontrivial invariant vectors then the nonsingular Gaussian dynamical system $(\Cal H,\mu_B,T_{(f,V)})$ is either conservative or totally dissipative. 
In fact, if $(Y,\goth C,\nu, S)$ is an ergodic conservative nonsingular dynamical system then the direct product
 $T_{(f,V)}\times S$ is either conservative or totally dissipative.
\endproclaim

The following theorem was first proved in \cite{ArIsMa} in the case of mixing $V$.
We extend it to the mildly mixing case with a different  proof.

\proclaim{Theorem 3.12} Let $T_{(0,V)}$ be mildly mixing and let $f$ not be a $V$-coboundary. If   
$T_{(f,V)}$ 
 is conservative then the Maharam extension of $T_{(f,V)}$ is sharply weak mixing.
 In particular,  $T_{(f,V)}$ is of type $III_1$.
 \endproclaim

\demo{Proof}
Since $T_{(f,V)}$ is conservative, it follows from \thetag{3-7} and Proposition~2.3 that $T_{(f,V)}$ is sharply weak mixing.
%The transformation $ T_{(f_r,V)}$ preserves invariant the $\mu_r$-equivalent probability measure
%$\nu_r$ such that 
Let $a_r,\mu_r$ and $\nu_r$ be as above in this subsection.
Since $\frac{d\mu_r}{d\nu_r}=\text{exp}_{-a_r}$ for each $r\in\Bbb N$,
it follows from \thetag{3-1} and \thetag{3-4} that the distribution $\psi_r$ of $\log\frac{d\mu_r}{d\nu_r}$ defined on  $(\Cal H_{0,r},\nu_r)$ is $\Cal N_{-\|a_r\|^2_0/2,\|a_r\|^2_0}$.
Hence, for all $m>n$, we have
$$
\align
\psi_{n+1}*\cdots*\psi_m&=\Cal N_{-0.5\sum_{r=n+1}^m\|a_r\|^2_0,\, \sum_{r=n+1}^m\|a_r\|^2_0}\quad\text{and}\\
\psi_{n+1}*\cdots*\psi_m*\delta_a&=\Cal N_{a-0.5\sum_{r=n+1}^m\|a_r\|^2_0,\, \sum_{r=n+1}^m\|a_r\|^2_0}
\endalign
$$
for each $a\in\Bbb R$.
We are going to show that the sequence  $(\psi_r)_{r=1}^\infty$  is ATQI.
First, it is straightforward to verify that for each $\sigma\in\Bbb R$ and $b>0$,
$$
\log\bigg(\frac{d\Cal N_{b-\sigma^2/2,\sigma^2}}{d\Cal N_{-\sigma^2/2,\sigma^2}}(t)\bigg)=
\frac {b(2t+\sigma^2-b)}{2\sigma^2}=\frac {bt}{\sigma^2}+
\frac{b(\sigma^2-b)}{2\sigma^2},\quad t\in\Bbb R.
$$
Hence, if  $t\ge-\sigma^2$ and $\sigma^2\ge 2b$ then
$$
\frac{d\Cal N_{b-\sigma^2/2,\sigma^2}}{d\Cal N_{-\sigma^2/2,\sigma^2}}(t)\ge e^{-b+\frac {b(\sigma^2-b)}{2\sigma^2}}>e^{-\frac{3b}4}.\tag3-8
$$
Moreover,
$$
\align
\int_{-\sigma^2}^{+\infty}{d\Cal N_{-\sigma^2/2,\sigma^2}}(t)&=\frac 1{\sigma\sqrt{2\pi}}
\int_{-\sigma^2}^{+\infty}
e^{-\frac12\Big(\frac{t+\sigma^2/2}{\sigma}\Big)^2}dt\\
&=
\frac 1{\sigma\sqrt{2\pi}}
\int_{-\sigma^2/2}^{+\infty}
e^{-\frac{t^2}{2\sigma^2}}dt\\
&=
\frac 1{\sqrt{2\pi}}
\int_{-\sigma/2}^{+\infty}
e^{-\frac{t^2}{2}}dt,
\endalign
$$
i.e. $\Cal N_{-\sigma^2/2,\sigma^2}((-\sigma^2,+\infty))=\Cal N_{0,1}((-\sigma/2,+\infty))$.
%Applying these estimates  to $\psi_{n+1},\dots,\psi_m$, we obtain the following:
Obviously, we have
$$
\psi_{n+1}*\cdots*\psi_m*\delta_a\sim \psi_{n+1}*\cdots*\psi_m\quad\text{for all $n<m$.}
$$
We now set $\zeta_a:=e^{-\frac{3a}4}$.
Next, we note that   $f$ is not a $V$-coboundary if and only if $\sum_{r=1}^\infty\|a_r\|_0^2=\infty$.
Hence for each $n>0$, there is $m>n$ such that $\sum_{r=n+1}^m\|a_r\|^2_0>2a$.
Let $W_{n+1,m}:=\big[ -\sum_{r=n+1}^m\|a_r\|^2_0,+\infty\big)\subset\Bbb R$.
Then \thetag{3-8} yields that
$$
\frac{d(\psi_{n+1}*\cdots*\psi_m*\delta_a)}{d(\psi_{n+1}*\cdots*\psi_m)}(t)\ge \zeta_a\quad\text{for all $t\in W_{n+1,m}$.}
$$
Moreover, $(\psi_{n+1}*\cdots*\psi_m)(W_{n+1,m})=\Cal N_{0,1}\Big(\Big(-0.5\sqrt{\sum_{r=n+1}^m\|a_r\|^2_0},+\infty\Big)\Big)\approx 1$ if $m$ is large.
Hence $(\psi_r)_{r=1}^\infty$  is ATQI.
It follows now from Theorem~2.10 that the Maharam extension of $T_{(f,V)}$ is sharply weak mixing.
\qed
\enddemo

\subhead 3.6. One-parametric  family of nonsingular Gaussian  systems
\endsubhead
We note that \thetag{3-5} determines a one-to-one homomorphism $\Bbb R^*\ni t\mapsto \alpha_t$ from the multiplicative group
$\Bbb R^*$ to the  group of continuous automorphisms of Aff$(\Cal H_0)$.
Therefore for each $A\in\text{Aff}(\Cal H_0)$, one can consider a one-parametric family
of nonsingular Gaussian transformations $T_{\alpha_t(A)}\in\text{Aut}(\Cal H,\mu_B)$, $t\in\Bbb R^*$\footnote{We note that the map  $\Bbb R^*\ni t\mapsto T_{\alpha_t(A)}$ is not a group homomorphism}.
Our purpose is this section is to investigate how the dynamical properties of $T_{\alpha_t(A)}$ depend on $t$.
It is straightforward to verify that the linear operator $-I$ of $\Cal H$ preserves $\mu_B$ and conjugates $T_{\alpha_t(A)}$ with $T_{\alpha_{-t}(A)}$.
Therefore it suffices to consider only the transformations $T_{\alpha_t(A)}$ with $t\in\Bbb R^*_+$.
\comment

Given $a\in\Cal H_0$ and $t>0$, it is straightforward to deduce from \thetag{3-3} that $(\exp_a)^t\in L^1(\Cal H,\mu_B)$,
$$
 \|(\exp_a)^t\|_1=e^{\|a\|_0^2\frac{t^2-t}2}\qquad\text{and}\qquad (\exp_a)^t/\|(\exp_a)^t\|_1=\exp_{ta},
$$
 where the 1-norm is taken with respect to $\mu_B$.
 If $f=a-Va$ then $\nu=\exp_{-a}\mu_B$ is $T_{(f, V)}$-invariant.
 And $\mu_B^{(t)}= \exp_{(t-1)a}\mu_B$ and $|\phi^t\|_1=e^{(t-1)\|a\|^2}$.

Now, if $(f,V)\in\text{Aff}(\Cal H_0)$ and $V$ has no non-trivial invariant vectors 
then we can represent $T_{(f,V)}$ as IDPFT (see~\thetag{3-7}).
There is an orthogonal  decomposition $\Cal H_0=\bigoplus_{r=1}^\infty\Cal H_{0,r}$ of $\Cal H_0$ into a sum of $V$-invariant subspaces $\Cal H_{0,r}$ such that
 $f=\bigoplus_{r=1}^\infty f_r$ for some $V$-coboundaries $f_r\in\Cal H_{0,r}$ for each $r>0$.
Since $\frac{d\mu_B\circ T_{(f,V)}^{-1}}{d\mu_B}=\exp_f$
$$
\prod_{n=1}^\infty\|(\exp_{-a_r})^t\|_1=e^{\frac{(t^2-t)}{2}\sum_{r=1}^\infty \|a_r\|_0^2}
=e^{\frac{(t^2-t)}{2} \|f\|_0^2}<+\infty,
$$ 
it follows that Cal$(\mu_B)=[1,+\infty)$ and $(\Cal H,\mu_B^{(t)}, T_{(f,V)})=(\Cal H,\mu_B, T_{(tf,V)})$.
\endcomment

\proclaim{Proposition 3.13 \cite{ArIsMa}} Given $A=(f,V)\in\text{\rom{Aff}}(\Cal H_0)$ such that $V$ has no non-zero invariant vectors, there is $t_{\text{\rom{diss}}}(A)\in [0,+\infty]$ such that the transformation $T_{\alpha_t(A)}$ is conservative if $0<t<t_{\text{\rom{diss}}}(A)$ and totally dissipative if $t>t_{\text{\rom{diss}}}(A)$.
\endproclaim
\demo{Proof} Let $A=(f, V)$ with $f\in\Cal H_0$ and $V\in\Cal O(\Cal H_0)$.
It is sufficient to show that if $T_A$ is totally dissipative then for each $t>1$, the Gaussian transformation $T_{\alpha_t(A)}$ is totally dissipative.
Since $T_A$ is totally dissipative, the Hopf criterion yields that
$
\sum_{n=0}^\infty\frac{d\mu_B\circ T_A^n}{d\mu_B}(x)=
\sum_{n=0}^\infty
e^{\langle B^{-1}f^{(n)},x\rangle-\frac12\|f^{(n)}\|^2_0}<\infty$
 for $\mu_B$-a.e. $x\in\Cal H$.
Hence,  there is $N_x>0$ such that $\langle B^{-1}f^{(n)},x\rangle-\frac12\|f^{(n)}\|^2_0<0$ for all $n>N_x$.
It follows that
$$
\langle tB^{-1}f^{(n)},x\rangle-\frac{t^2\|f^{(n)}\|^2_0}2<t\bigg(\langle B^{-1}f^{(n)},x\rangle-\frac{\|f^{(n)}\|^2_0}2\bigg)< \langle B^{-1}f^{(n)},x\rangle-\frac{\|f^{(n)}\|^2_0}2
$$
for all $n>N_x$.
Hence $
\sum_{n=0}^\infty
e^{\langle tB^{-1}f^{(n)},x\rangle-\frac12\|tf^{(n)}\|^2_0}<\infty$
 for $\mu_B$-a.e. $x\in\Cal H$.
 Since $tf^{(n)}=(tf)^{(n)}$, we deduce from the Hopf criterion that $T_{\alpha_t(A)}$
 is dissipative, as desired. \qed
 \enddemo

We recall that the {\it Poincare exponent} of $A=(f,V)\in\text{Aff}(\Cal H_0)$ \cite{ArIsMa} is
$$
\delta_A:=\inf\bigg\{\alpha>0\mid\sum_{n=1}^\infty e^{-\alpha\|f^{(n)}\|^2_0}<+\infty\bigg\}\in[0,+\infty].
$$

%It was also shown there that if $0<t_1<t_2<t_{\text{diss}}(A)<+\infty$ then $T_{\alpha_{t_1}(A)}$ is not isomorphic to $T_{\alpha_{t_2}(A)}$.

For completeness of our argument we give a proof of the following proposition.

\proclaim{Proposition 3.14 \cite{ArIsMa}}
$\sqrt{2\delta_A}\le t_{\text{diss}}(A)\le 2\sqrt{2\delta_A}$.
\endproclaim
\demo{Proof \cite{ArIsMa}}
Let $t>t_{\text{diss}}(A)$.
Since $T_{\alpha_t(A)}$ is isomorphic to $T_{\alpha_{-t}(A)}$, the two transformations are dissipative.
Therefore, by the Hopf criterion, 
$$
\sum_{n=0}^\infty
e^{t\langle B^{-1}f^{(n)},x\rangle-\frac{t^2}2\|f^{(n)}\|^2_0}<\infty\quad\text{and}\quad
 \sum_{n=0}^\infty
e^{-t\langle B^{-1}f^{(n)},x\rangle-\frac{t^2}2\|f^{(n)}\|^2_0}<\infty
$$
at a.e. $x$. 
Since $e^{t\langle B^{-1}f^{(n)},x\rangle}+e^{-t\langle B^{-1}f^{(n)},x\rangle}\ge 2$ for each $x\in X$, it follows that $\sum_{n=0}^\infty e^{-\frac {t^2}2\|f^{(n)}\|^2_0}<\infty$, i.e. $\delta_A\le t^2/2$ and hence  $\delta_A\le t_{\text{diss}}(A)^2/2$.

On the other hand, if $t<t_{\text{diss}}(A)$ then $T_{\alpha_t(A)}$ is conservative and hence
$$
\sum_{n=0}^\infty
e^{t\langle B^{-1}f^{(n)},x\rangle-\frac{t^2}2\|f^{(n)}\|^2_0}=+\infty.
$$
Therefore, 
$$
+\infty=\sum_{n=0}^\infty\int_{\Cal H}e^{\frac 12t\langle B^{-1}f^{(n)},x\rangle-\frac{ t^2}4\|f^{(n)}\|^2_0}d\mu_B(x)=
\sum_{n=0}^\infty e^{-\frac{t^2}8\|f^{(n)}\|^2_0}.
$$
Hence $\delta_A\ge\frac {t^2}8$ and therefore  $\delta_A\ge\frac {t_{\text{diss}}(A)^2}8$.
\enddemo

\comment
\subhead 3.7. Mariusz question
\endsubhead
Let $T\in
\text{Aut}(X,\mu)$ and let $S$ be an ergodic finite  measure preserving transformation of a standard probability space $(Y,\nu)$.
We say  that a function $f\in L^\infty(X,\mu)\setminus\{0\}$
is an $L^\infty$-eigenfunction of $T$ with the eigenvalue $\lambda\in\Bbb T$ if $f\circ T=\lambda f$.
Denote by $e(T)$ the set of all eigenvalues of $T$.
Then $e(T)$ is a Borel subgroup of $\Bbb T$.

\proclaim{Proposition} Let $\sigma_S$ denote the reduced maximal spectral type of $S$.
If $\sigma_S(e(T))=0$ then every $(T\times S)$-invariant function $F\in L^\infty(X\times Y,\mu\otimes\nu)$ does not depend on the first coordinate.
\endproclaim

\demo{Proof} For $x\in X$ and $y\in Y$, we let $F_x(y):=F(x,y)$.
Then $F_x\in L^\infty(Y,\nu)\subset L^2(Y,\nu)$ for a.e. $x\in X$.
Without loss of generality we may assume that $\int_YF_x\,d\nu=0$,  i.e. $F_x\in L^2_0(Y,\nu)$ for a.e. $x\in X$.
Since $F\circ (T\times S)=F$, it follows  that
$F_{Tx}=U_SF_{x}$ for a.e. $x\in X$, where $U_S$ denote the Koopman operator generated by $S$.
By the spectral decomposition theorem, there is a unitary operator $V:L^2_0(Y,\nu)\to \int_\Bbb T^\otimes\Cal H_zd\sigma_S(z)$ denote 
r such that $VU_SV^*=\int^\otimes_{\Bbb T} zI_zd\sigma_S(z)$.
Consider a function 
$$
G:X\times\Bbb T\ni (x,z)\mapsto G(x,z):=VF_x(z)\in\Bbb C.
$$
Then $G(Tx,z)=zG(x,z)$  and hence
$|G(Tx,z)|=|G(x,z)|$ for  $\mu\otimes\sigma_S$-a.e. $(x,z)\in X\times\Bbb T$.
 We now set
 $$
 \widetilde G(x,z):=
 \cases
 \frac{ G(x,z)}{| G(x,z)|} &\text{ if }G(x,z)\ne 0,\\
 G(x,z) &\text{otherwise}.
 \endcases
 $$
 Then $ (x,z)\mapsto\widetilde G(x,z)$ is bounded map and $\widetilde G(Tx,z)=z\widetilde G(x,z)$ almost everywhere.
Hence the mapping $X\ni x\mapsto \widetilde G(x,z)$ is either $0$ or an $L^\infty$-eigenfunction of $T$ for each $z\in\Bbb T$.
The latter situation can happen only for a $\sigma_S$-null subset of $\Bbb T$ by the assumption of Proposition.
Therefore $\widetilde G(x,z)=0$ a.e. and hence 
 $\widetilde G(x,z)=0$ a.e.
 This, in turn yields that $F=0$, as desired.
 \qed
 \enddemo

\endcomment

\enddocument